\documentclass[10pt]{article}
\usepackage{a4,amssymb,verbatim,amsmath}

\usepackage{rotating}

\newtheorem{theorem}{Theorem}[section]
\newtheorem{lemma}[theorem]{Lemma}
\newtheorem{remark}[theorem]{Remark}
\newtheorem{corollary}[theorem]{Corollary}
\newtheorem{proposition}[theorem]{Proposition}

\newtheorem{ex}{Example}[section]
\newenvironment{example}{\begin{ex}\rm}{ \hfill $\Diamond$ \end{ex}
        \vskip4pt}
\newtheorem{ass}{Assumption}[section]

\numberwithin{equation}{section}

\begin{document}

\def\sso#1{\ensuremath{\mathfrak{#1}}}  % to make gotic letters
\def\ssobig#1{{{#1}}}

\newcommand{\ddx}{\partial \over \partial x}
\newcommand{\ddy}{\partial \over \partial y}

\newcommand{\dx}{\partial  _x}
\newcommand{\dy}{\partial  _y }

\begin{center}
{\bf \Large    Second-order delay ordinary differential equations, their symmetries and application to a traffic problem}
\end{center}

\bigskip

\begin{center}
{\large Vladimir A. Dorodnitsyn}$^a$,
{\large Roman Kozlov}$^b$,  \\
{\large Sergey V. Meleshko}$^c$
{\large and
Pavel Winternitz}$^d$
\end{center}

\bigskip

\noindent $^a$
Keldysh Institute of Applied Mathematics, Russian Academy of Science, \\
Miusskaya Pl.~4, Moscow, 125047, Russia; \\
{e-mail:  Dorodnitsyn@Keldysh.ru } \\
%  .... \\[10pt]
$^b$
Department of Business and Management Science, Norwegian School of Economics,
Helleveien 30, 5045, Bergen, Norway;  \\
{e-mail: Roman.Kozlov@nhh.no}  \\
%  .... \\[10pt]
$^c$
School of Mathematics, Institute of Science,
Suranaree University of Technology, 30000, Thailand; \\
{e-mail: sergey@math.sut.ac.th} \\
%  .... \\[10pt]
$^d$
Centre de Recherches Math\'ematiques
and
D\'epartement de math\'ematiques et de statistique,
Universit\'e de Montr\'eal,
Montr\'eal, QC, H3C 3J7, Canada; \\
{e-mail: wintern@crm.umontreal.ca} \\

\begin{center}
{\bf Abstract}
\end{center}
\begin{quotation}
This article is the third in a series the aim of which is to
use Lie group theory to obtain exact analytic solutions of
Delay Ordinary Differential Systems  (DODSs).
Such a system consists of
two equations involving one  independent variable $x$ and one
dependent variable $y$. As opposed to ODEs the variable $x$ figures
in more than one point (we consider the case of two points, $x$ and
$x_-$). The dependent variable $y$ and its derivatives figure in
both $x$ and $x_-$.
Two previous articles were devoted to {\it first}-order DODSs,
here we concentrate on a large class of {\it  second}-order ones.
We show that within this class the symmetry
algebra can be of dimension $n$ with  $0 \leq n \leq 6$ 
for nonlinear DODSs and must be $n=\infty$ for linear or linearizable
ones. 
The symmetry algebras can be used to obtain exact particular group invariant solutions. 
As a
specific application we present some exact solutions of a
DODS model of traffic flow.
% A group classification of second-order delay ordinary differential
% equations (DODEs) accompanied by an equation for delay parameter
% (delay relation) is presented. A subset of such systems
% (delay ordinary differential systems or DODSs)
% which consists of linear DODEs and solution-independent delay
% relations have infinite-dimensional symmetry algebras, as do
% nonlinear ones that are linearizable by an invertible transformation
% of variables.
% Genuinely nonlinear DODSs have symmetry algebras of
% dimension $n$, $0 \leq n \leq 6$.
% {\bf It is shown how exact analytical solutions of invariant DODSs
% can be obtained using symmetry reduction.   }
\end{quotation}

% \bigskip

% \noindent {08.07.2020}

% \eject

\section{Introduction and formulation of the problem}

\label{Introduction}

Two previous articles were devoted to the adaptation  of  Lie group and Lie algebra theory
to the study  of delay ordinary differential equations~\cite{DorKozMelWin2018a, DorKozMelWin2018b}.
 In these articles we restricted ourselves to the case of first-order DODEs, supplemented by a general delay equation.
Thus we considered first-order delay ordinary differential systems (DODSs) of the form
\begin{subequations}     \label{first_order_DODS}
\begin{gather}  
\label{first_order}
\dot{y} = f ( x, y, y_-  ),  
\qquad  
{ \partial f  \over  \partial y_- }  {\not\equiv}  0  , 
\qquad 
x \in I ,  
\\
x_- = g ( x, y, y_-   ) ,
\qquad 
x_- <  x , 
\qquad 
g  ( x, y, y_- ) {\not \equiv}
\mbox{const} .
\end{gather}  
\end{subequations}  
Here  $I$ is a finite or {semi}finite interval and $f$ and $g$ are arbitrary smooth functions.
For details and motivation we refer to our articles~\cite{DorKozMelWin2018a, DorKozMelWin2018b}.

Our main results were the following:

\begin{enumerate}

\item

We  classified DODSs of the form~(\ref{first_order_DODS})  into conjugacy classes
under arbitrary Lie point transformations and found that their Lie point symmetry groups
can have dimension $ n$ with $ 0 \leq n \leq 3$,   or infinity.
If  $n=\infty$ the DODE is linear or linearizable by a point transformation.
In general the Lie algebra of the infinite-dimensional symmetry group of a linear DODS is a solvable Lie algebra
with an infinite-dimensional {nil}radical and is realized by vector fields of the form
\begin{equation*}
 X= \eta (x,y)    {\ddy}   .
\end{equation*}

\item

If the symmetry algebra of a DODS contains a 2-dimensional {sub}algebra
realized by linearly connected vector fields,  then this DODS  is linearizable (or already linear) with $g=g(x)$.

\item

The symmetry algebra  for genuinely nonlinear DODEs  has  dimension $n  \leq 3 $.
For algebras with $n=2$ or $3$  we presented a method for obtaining exact particular solutions.

\end{enumerate}

% The purpose of this article is to perform a similar analysis for second-order DODSs
% and to show how the previously obtained results can be generalized to more general DODSs.

% In the recent article~\cite{DorKozMelWin2018a}  there was presented
% the Lie group classification of first-order delay ordinary differential equations.
% We refer to~\cite{DorKozMelWin2018a}  and the references therein
% for motivation and a brief survey of the field.
% Linear  first-order delay ordinary differential equations
% were treated in detail in a separate article~\cite{DorKozMelWin2018b}.

% \section{OLD VERSION}
%  (It is modified text from our first paper)}

% \label{Formulation_of_the_problem}

Here we provide a similar analysis of
SECOND-ORDER  delay ordinary differential systems

\begin{subequations}   \label{DODS}
\begin{gather}  
  \label{DODE}
\ddot{y} = f ( x, y, y_- , \dot{y} , \dot{y} _-   ),
\qquad
\left( { \partial f  \over  \partial y_- } \right) ^2    +   \left( { \partial f  \over  \partial \dot{y}_- } \right) ^2
{\not\equiv}  0  ,
\qquad
x \in I , 
\\
% where $I \subset  \mathbb{R} $ is some finite or semi{infinite} interval.
% We will be interested in symmetry properties of this DODE,
% which is considered locally, independently of initial conditions.
% For equation~(\ref{DODE})  we have to specify the delayed point $x_-$
% where the delayed function value $ y_- = y(x_-) $
% and the derivative value $ \dot{y}_- = \dot{y} (x_-) $ are taken,
% otherwise the problem is not fully determined.
% To determine the problem completely  we supplement the DODE with a delay relation
    \label{delay}
x_- = g ( x, y, y_- , \dot{y} , \dot{y} _-   ) ,
\qquad
x_- <  x ,
\qquad
g  ( x, y, y_- , \dot{y} , \dot{y} _-  ) {\not  \equiv}   \mbox{const} . 
\end{gather}
\end{subequations}
Thus a DODS consists of a delay  ordinary differential equation (DODE~(\ref{DODE}))  
and a delay equation~(\ref{delay})

% We will use the same terminology as in our previous two articles,
% calling the above equations a DODS,
% i.e. a differential ordinary delay system.
% For convenience we sometimes use  the notation
% $ \Delta x = x - x_- $ and reformulate~(\ref{delay}) in terms of $ \Delta x $.

In most of the existing literature the delay $  \Delta x =  x - x_- = \tau >  0 $  is constant.
An alternative is to impose specific conditions characterizing the physical problem under consideration.
We impose the conditions given in~(\ref{delay}) and leave the function $g$  free
and to be determined by symmetry considerations.

As in the case or ordinary differential equations  (ODEs) we will be working with vector fields
of the form
\begin{equation}    \label{operator1}
X _{\alpha}  =   \xi  _{\alpha}   (x,y)  {\ddx}
+ \eta  _{\alpha}  (x,y)  { \ddy}  .
\end{equation}
Integrating them,  we obtain (local) Lie point transformations
acting on the independent variable $x$, the dependent one $y$
and thus also on functions $ y(x) $.
They must be prolonged to act on DODSs viewed as functions of
$x$, $y$, $\dot{y}$ and $\ddot{y}$
evaluated at two points, the reference point $x$ and the delay point $x_-$.
The transformations will leave the DODS~(\ref{DODS})  
invariant on its solution set, i.e. the prolongation of the vector field~(\ref{operator1})
must annihilate the DODS on its solution set.
From this point of view the symmetry group of DODS are reminiscent
of those of ordinary difference systems that leave invariant the solution set of two equations:
the difference equation itself and an equation determining
the lattice \cite{Dorodnitsyn2011, DKW2000, DKW2004,  Levi2000, Levi2001, Levi1997, Levi1991, Levi2006}.

Delay differential equations have many properties that 
distinguish them from differential equations without delay 
\cite{Driver1977, bk:Elsgolts[1955],   Hale1977, Kolmanovskii, bk:Wu1996}. 
For earlier work on solutions of delay equations 
be they numerical,  or exact, 
using group theory 
and other approaches  we refer to 
\cite{Grigoriev2010,  LongMeleshko2017, Nass2019a,  Nass2019b,  Ozbenli2020, 
bk:PolyaninZhurov[2014]a, bk:PolyaninZhurov[2014]b, PolyaninSorikin, PolyaninZhurov2015, 
Smith2011,  Tanthanuch2012}.  
We recall that DODEs are standardly solved numerically or otherwise using
the {\it method of steps} as described e.g. in~\cite{Myshkis1989}
and adapted to {non}constant delay in our previous articles~\cite{DorKozMelWin2018a, DorKozMelWin2018b}.

In Section~\ref{general_theory} we describe our general classification procedure.
All nonlinear DODSs of the type~(\ref{DODS})   are classified
into symmetry classes in  Section~\ref{Classification}.
In Section~\ref{linearly_connected_section} we investigate linear DODSs.
An algorithm for calculating exact analytic solutions of DODSs is given in  Section~\ref{Invariant solutions}.
Applications to a "follow the leader" model of traffic flow are given in Section~\ref{micro_model}.
Finally, the concluding remarks are presented in
Section~\ref{Conclusion}.

\section{Second-order DODSs invariant under local Lie point transformations groups}

%  (It is modified text from our first paper ???)

\label{general_theory}

The method of constructing invariant first-order DODSs
in Refs.~\cite{DorKozMelWin2018a, DorKozMelWin2018b}
can be generalized to DODSs of any order and any number of delay points.
For the DODS~(\ref{DODS})  we prolong the vector field~(\ref{operator1}) to
\begin{equation}    \label{prolongation}
\mbox{\bf pr} X  _{\alpha}
=    \xi _{\alpha}  {\ddx}
+ \eta _{\alpha}   { \ddy}
+   \xi  _{\alpha}   ^-   {\partial  \over \partial x_-}
+ \eta _{\alpha}   ^-  {\partial   \over \partial y_-}
+   \zeta  _{\alpha} ^1  {\partial  \over \partial \dot{y} }
+   \zeta  _{\alpha} ^{1-}  {\partial  \over \partial \dot{y}_- }
+  {\zeta}  _{\alpha}  ^2   {\partial  \over \partial \ddot{y} }
\end{equation}
with
\begin{equation*}
\xi _{\alpha}   = \xi  _{\alpha} (x,y) ,
\qquad
\eta _{\alpha}  = \eta   _{\alpha}  (x ,y )  ,
\qquad
\xi _{\alpha} ^-    = \xi   _{\alpha}  (x_- ,y_-)  ,
\qquad
\eta  _{\alpha} ^-    = \eta   _{\alpha} (x_- ,y_-)  ,
\end{equation*}
\begin{equation*}
\zeta   _{\alpha} ^1  (x,y, \dot{y} )  = D  ( \eta  _{\alpha}  )  - \dot{y}   D  ( \xi   _{\alpha} )  ,
\qquad
\zeta   _{\alpha}  ^{1-} (x_-,y_-, \dot{y}_- )  = D  _- ( \eta  _{\alpha}  ^-  )  - \dot{y}_-    D _-  ( \xi   _{\alpha} ^- )  ,
\end{equation*}
\begin{equation*}
{\zeta}   _{\alpha}  ^2   (x,y, \dot{y},  \ddot{y})  = D ( \zeta  _{\alpha} ^1  )  - \ddot{y}   D  ( \xi   _{\alpha} )  ,
\end{equation*}
where $D$ and $D_-$ are the total derivative operators in points $ x $ and $x_-$, respectively. 
Simply put,  the coefficients   $ \zeta  _{\alpha} ^1  $ and  $   {\zeta}  _{\alpha}  ^2  $
are calculated as for ODEs~\cite{Olver1986, Ovsiannikov1982} 
and the coefficients  $ \xi  _{\alpha}   ^-   $,  $  \eta _{\alpha}   ^-  $ and $   \zeta  _{\alpha} ^{1-}  $
are obtained by  shifting $ x$  and $ y $ to   $x_- $ and  $ y_- $
in   $ \xi  _{\alpha}     $,  $  \eta _{\alpha}    $ and $   \zeta  _{\alpha} ^{1}  $.

We construct, classify and represent all symmetry classes of DODSs  of the type~(\ref{DODS})
using essentially the same method as we used in previous articles~\cite{DorKozMelWin2018a, DorKozMelWin2018b}
for first-order DODSs:

\begin{enumerate}

\item

For each algebra of vector fields in the list given in~\cite{Gonzalez1992}
we construct the prolongations~(\ref{prolongation})
of the chosen basis vectors $ \{ X_{\alpha},  \alpha = 1, ..., n \} $.

\item

Running through all the algebras of the list,
we construct the "strongly invariant" DODSs out of universal point invariants satisfying
\begin{equation}    \label{strong}
\mbox{\bf pr} X  _{\alpha}    \Phi ( x, y ,x_- , y_- ,\dot{y} ,   \dot{y}_-  ,  \ddot{y} ) = 0 ,
\qquad
 \alpha = 1, ... , n  .
\end{equation}
The method of characteristics gives us a set of independent elementary invariants
of the corresponding Lie group action.
We label them $ I_1 $, ... $ I_k $, where $ k $ satisfies
\begin{equation}    \label{dimension_m}
k  =  \mbox{dim} \  M - (   \mbox{dim} \   G  -    \mbox{dim} \   G  _0 )
\end{equation}
with   $  M  \sim ( x, y ,x_- , y_- ,\dot{y} ,  \dot{y}_- ,  \ddot{y} )  $.
In~(\ref{dimension_m})    $ G $  is the local Lie point symmetry group
corresponding to the considered Lie algebra
and  $G _0  \subset G $ is the stabilizer of a generic point on the manifold $M$.

The strongly invariant DODSs can all be written as

\begin{equation}   \label{implicit1}
F (   I_1, ... , I_k ) = 0   ,
\qquad
G (   I_1, ... , I_k ) = 0   ,
\qquad
\mbox{det}  \left(   {\partial ( F , G )  \over \partial ( \ddot{y}  , x_- ) }  \right)  \neq 0
 %  \label{implicit2}
\end{equation}
($F$  and  $ G$  are otherwise arbitrary smooth functions).

It is convenient to introduce a matrix $ Z $  for each algebra of prolonged vector fields,
namely
\begin{equation}
Z = \left(
\begin{array}{ccccccc}
\xi_1    & \eta_1   &  \xi _1 ^-  &  \eta  _1 ^-   & \zeta _1 ^1  &   \zeta _1 ^{1-}    &  {\zeta}   _1   ^2 \\
\vdots  & & & &  \\
\xi_n    & \eta_n   &  \xi _n ^-  &  \eta  _n ^-   & \zeta _n  ^1 &   \zeta _n ^{1-}    &  {\zeta}   _n   ^2  \\
\end{array}
\right)   .
\end{equation}
The number of strong invariants in then expressed as
\begin{equation}
k =   \mbox{dim} \   M -  \mbox{rank} \  Z
% \qquad
% k   \geq  0   ,
\end{equation}
where   $ \mbox{rank} \  Z  $ is evaluated at a generic point.

\item

Complement  the strongly invariant DODSs
found in Step~2 above by weakly invariant ones.
These are found using "weak invariants".
They satisfy a weaker equation than~(\ref{strong}), namely
 \begin{equation}    \label{weak}
\left.
\mbox{\bf pr} X  _{\alpha}    \Phi ( x, y ,x_- , y_- ,\dot{y} ,   \dot{y}_-  ,  \ddot{y} )
\right |  _{  \Phi = 0 } = 0 ,
\qquad
 \alpha = 1, ... , n  .
\end{equation}
 The weak invariants   $ J_a $ are then invariant for a Lie group action on a {sub}manifold of $ M $:
$ G ( J_a ) \neq J_a $, but  $ \left. G ( J_a ) \right | _{ J_a  = 0 } = 0 $.

\end{enumerate}

\section{Representative list of second-order  {non}linear  DODSs}

\label{Classification}

We proceed as outlined in Section~\ref{general_theory}
and as we did in our previous article on first-order DODSs~\cite{DorKozMelWin2018a}.
Thus we run through the list of Lie algebras of vector fields presented in~\cite{Gonzalez1992},
proceeding by dimension.
For each algebra in the list we construct a basis for all strong and weak invariants in the space
with coordinates  $ ( x, y ,x_- , y_- ,\dot{y} ,  \dot{y}_- ,  \ddot{y} )  $.
We then construct the most general DODS of the form~(\ref{DODS}).

We shall not present any of the standard calculations here.
The final results are summed up in Tables 1, 2, 3 and 4.
The algebras have dimensions  $ \mbox{dim} \  L $ satisfying     $ 1 \leq \mbox{dim} \  L  \leq 6 $.
The reason for this cut-off is that only linear and linearizable second-order DODSs
(of the form~(\ref{DODS})) have symmetry algebras of dimensions larger than  $  \mbox{dim} \  L  = 6 $.

Linear DODSs are studied below in Section~\ref{linearly_connected_section}.
Contrary to second-order linear ODEs (with no delay) linear DODSs of any order allow
infinite-dimensional Abelian Lie point symmetry algebras,
realized by linearly connected vector fields.
These algebras are just expressions of the linear superposition principle for linear DODSs.

Let us now discuss the individual tables.
Table~1 corresponds to algebras of dimension 1 and 2.
In general the delay equation as well as the DODE is nonlinear.
For special choices of the functions $f$ and $g$ they may be linear
and then the symmetry algebra is infinite-dimensional.
Notice that for     $\mbox{\bf A}_{2,1}$ and $\mbox{\bf A}_{2,3}$
the vector fields are linearly connected.
Contrary to the case of first-order DODSs this does not imply that
the second-order DODSs are linear.

From Table~2 and Ref.~\cite{SWbook} we see that each isomorphy class
of three-dimensional real Lie algebras allows at least one class of invariant DODSs
and the DODEs and delay equations involve arbitrary functions $f$ and $g$ of two variables.
We mention that  
the particular case $ a=1$ of algebra   $\mbox{\bf A}_{3,3} ^a$,    
one of the algebras
 $ {\sso {sl}}  (2, \mathbb{R})   $ 
(namely  $\mbox{\bf A}_{3,11}  $), 
and the algebra   $ 3  {\sso {n}}_{1,1}  $ 
are realized entirely by linearly connected vector fields 
since all three vector fields involve only $y$-derivatives. 
This does not force the DODSs to be linear.

The classification results for  $  \mbox{dim} \  L  = 4 $ are summed up in Table~3.
{All}together 13 classes of four-dimensional 
 indecomposable real Lie algebra exist~\cite{SWbook}. 
Two of them $ {\sso {s}}_{4,7}  $ and $ {\sso {s}}_{4,9}  $ are missing in Table~3 
because they cannot be represented by real vector fields in two dimensions~\cite{Gonzalez1992, Nesterenko2}. 
The  decomposable algebras 
$ {\sso {n}}_{1,1} \oplus {\sso {n}}_{3,1} $
and
$ {\sso {n}}_{1,1} \oplus {\sso {o}}  ( 3 )  $ 
are missing for the same reason. 
Two other Lie algebras are omitted in Table~3 for a different reason. 
They are realized by linearly connected vector fields 
and we shall see in Section~\ref{linearly_connected_section} 
that they lead to linear (or linearizable) DODSs.  
They correspond to the solvable Lie algebra $ {\sso {s}}_{4,3} $
with basis $ \{  \partial _y  , x \partial _y , 
%   \} \oplus \{  
\chi (x)    \partial _y  , y \partial _y   \}  $, 
$ \ddot{\chi} (x)  {\not \equiv}  0 $ 
and the Abelian algebra  
$ 4 {\sso {n}}_{1,1}  $ 
with basis 
$ \{  \partial _y  , x \partial _y , \chi _1 (x)    \partial _y  , \chi _2 (x)   \partial _y   \}  $, 
where the functions  $ \{   1  , x  , \chi _1 (x)   , \chi _2 (x)    \}  $ 
are linearly independent.

% The classification results for  $  \mbox{dim} \  L  = 4 $ are summed up in Table~3.
% All isomorphism classes of four-dimensional Lie algebras are represented except
% the indecomposable solvable Lie algebra
% $\mbox{\bf A}_{4,7}  $ and $\mbox{\bf A}_{4,9} $
% { \bf  (are indices correct ????) }
% with nonzero commutation relations
% $$
% [  e_2 , e_3  ] = e_1 ,
% \qquad
%  [  e_4 , e_2  ] = - e_3 ,
% \qquad
%  [  e_4 , e_3  ] = e_2
% $$
% and
% $$
% [  e_2 , e_3  ] = 2 e_1 ,
% \qquad
%  [  e_4 , e_2  ] = \alpha e_2 - e_3 ,
% \qquad
%  [  e_4 , e_3  ] = e_2 +   e_3
% $$
% respectively.
% These are the only two indecomposable Lie algebras
% that can be realized by real vector fields in two dimensions.
% The decomposable Lie algebras missing in Table~3
% that cannot be realized are
% $ {\sso {n}}_{1,1} \oplus {\sso {n}}_{3,1} $
% and
% $ {\sso {n}}_{1,1} \oplus {\sso {O}}  ( 3 )  $.
% Two algebras that can be realized are $ 4 {\sso {n}}_{1,1}  $
% $ \{  \partial _y  , x \partial _y , \chi _1 (x)    \partial _y  , \chi _2 (x)   \partial _y   \}  $
% and
% $ {\sso {s}}_{2,1} \oplus 2 {\sso {n}}_{1,1} $
% $ \{  \partial _y  , x \partial _y   \} \oplus \{  \chi _1 (x)    \partial _y  , \chi _2 (x)   \partial _y   \}  $
% { \bf  (simplification needed  ????)}.
% The do not figure in Table~3 because they only lead to linear DODSs.

% and the decomposable nilpotent Lie algebra
% $   \{  e_1 , e_2 , e_4 \} + \{   e_4  \} $  (with $[  e_2 , e_3  ] = e_1$).
% These two algebras cannot be realized by vector fields in two dimensions.

All algebras in Table~3 correspond to nonlinear DODEs and delay relations
that both involve arbitrary functions $f$  and $g$ of one variable each.
For $f$  and $g$ generic the equations are nonlinear.
None of the symmetry algebras consists entirely of linearly connected vector fields,
though many of them have three-dimensional {sub}algebras  of linearly connected vector fields.

Finally, let us discuss Table~4 devoted to nonlinear DODSs with symmetry algebras
of dimension 5 and 6.
We present a representative list of all five- and six-dimensional Lie algebras
that can be realized by vector fields in two dimensions, 
do not contain a four-dimensional {sub}algebra   realized by linearly connected vector fields
and  do leave a DODS invariant.
For $  \mbox{dim} \  L  = 5 $ there are 8 such algebras.
Each class of DODSs depends on two arbitrary constants
but no arbitrary functions appear.
The Levi decomposable Lie algebra   $ {\sso {sl}}  (2, \mathbb{R})  \ltimes  2  {\sso {n}}_{1,1}  $
appears in two inequivalent ways.
The algebras $\mbox{\bf A}_{5,1}  $ -- $\mbox{\bf A}_{5,4}  $
are solvable with {nil}radical   $  {\sso {n}}_{4,1} $, 
i.e. $  \{  X_1 , ... ,  X_4 \}  $;
 $\mbox{\bf A}_{5,5}  $ is also solvable with
 {nil}radical   $  {\sso {n}}_{3,1} $, 
i.e. $  \{  X_1 , X_2 ,  X_3 \}  $;
 $\mbox{\bf A}_{5,6}  $ and $\mbox{\bf A}_{5,7}  $
 are isomorphic but not conjugate.
 The radical   $  \{  X_4 , X_5 \}  $
 is linearly connected for $\mbox{\bf A}_{5,6}  $
 but not for  $\mbox{\bf A}_{5,7}  $.

% {\bf NOT CORRECT:   all other isomorphy classes
% are indecomposable and solvable}.
% Each appears only once.

For $  \mbox{dim} \  L  = 6 $ there exist only three Lie algebras that fulfil all conditions.
The DODSs depend on one constant, no arbitrary functions. 
We mention that in cases   $\mbox{\bf A}_{5,1}  - \mbox{\bf A}_{5,5}  $  
and  $ \mbox{\bf A}_{6,2}  $
the DODE is actually linear, but the delay relation depends on the solution $y$
and this induces nonlinearity.

For $  \mbox{dim} \  L  \geq  7 $ the only Lie algebra realized
by the considered vector fields
that does not contain a four-dimensional {sub}algebra of linearly connected vector fields is
$ {\sso {sl}}  (3, \mathbb{R}) $ (see~\cite{Gonzalez1992}).
There is no second-order DODS invariant under the Lie group $ { {SL}}  (3, \mathbb{R}) $. 
All other Lie algebras with $  \mbox{dim} \  L  \geq  7 $
are contained in the series 20, ..., 28 of~\cite{Gonzalez1992}
and lead to linear DODSs.

\section{Symmetry algebras of linear second-order DODSs}

\label{linearly_connected_section}

In our previous articles~\cite{DorKozMelWin2018a, DorKozMelWin2018b}
on DODSs we considered in particular
the symmetries of linear and linearizable first-order DODSs.
It was shown that a special role is played by symmetry algebras of
linearly connected vector fields.
These are vector fields that in the case of a two-dimensional real manifold
can all be simultaneously reduced to the form
$$
X _{\alpha}  = \eta _{\alpha}  (x,y) {\ddy} .
$$

The representative list of algebras of vector fields in a real plane
(\cite{Gonzalez1992}, Table~1)
contains 3 individual algebras of this type, namely algebras 9,  10  and  11  in the table 
(here we use notations of~\cite{Gonzalez1992}),
and two infinite families,  20  and  21  in the same table.

We adapt these results to our needs. We include algebras
9  and   10  in the families  20  and  21, respectively.
The two families  now are
\begin{equation}   \label{series_S}
S_m =  \{  \partial _y  , x \partial _y  ,
\chi _2 (x)    \partial _y  , ..., \chi _m (x)   \partial _y   \}  ,
\qquad
0 \leq m < \infty   , 
\end{equation}
\begin{equation}    \label{series_H}
H_m =  \{  \partial _y  , x \partial _y  ,
\chi _2 (x)    \partial _y  , ..., \chi _m (x)   \partial _y ,  y \partial _y  \}  ,
\qquad
0 \leq m < \infty   , 
\end{equation}
where we put
$  \chi _0 (x) = 1 $, $  \chi _1 (x) = x $
and the functions $  \chi _0 (x) $, ..., $  \chi _m (x)  $
are linearly independent.   
We remark that algebra 11, namely  
 $ {\sso {sl}} ( 2, \mathbb{R}  ) \sim    \{  \partial _y  , y \partial _y  ,  y^2    \partial _y  \} $, 
is not included in the infinite families.

It was shown in Ref.~\cite{DorKozMelWin2018a} that
if a symmetry algebra of a first-order DODE contains
a two-dimensional {sub}algebra of linearly connected vector fields,
it  is linearizable.
That is the case of
$ S_1 =  \{  \partial _y  , x \partial _y    \}   $
for linear inhomogeneous DODSs or
$ H_0 =  \{  \partial _y  ,   y \partial _y  \} $
for homogeneous ones.
Here we show that for second-order DODSs (with one delay point)
a four-dimensional {sub}algebra of linearly connected vector fields
has the same effect.

Thus, we now turn to the general second-order DODS
with one delay point $ x_- $
\begin{equation}    \label{inhomogeneous}
 \ddot{y} =  a_1  ( x  )  \dot{y}  +   a_2   ( x  )  \dot{y} _-  +  a_3   ( x  )  y +  a_4   ( x  )  y_-  +  b  ( x  )  ,
\qquad
x_-   = g  ( x )   ,
\end{equation}
where $ a_i (x) $, $b(x)$  and  $ g(x) $ are arbitrary real functions,
smooth in some interval $ x \in I $. 
The functions  satisfy 
\begin{equation}
a_2^2    (x) + a_4^2   (x)   {\not\equiv} 0 ,
\qquad
g (x)  <  x,
\qquad
g (x)  {\not\equiv}  \mbox{const}  .
\end{equation}
We will mainly use homogeneous DODS, i.e.~(\ref{inhomogeneous}) with $ b  ( x  ) = 0 $  
\begin{equation}    \label{homogeneous}
 \ddot{y} =  a_1  ( x  )  \dot{y}  +   a_2   ( x  )  \dot{y} _-  +  a_3   ( x  )  y +  a_4   ( x  )  y_-    ,
\qquad
x_-   = g  ( x )   . 
\end{equation}

\begin{theorem}  \label{linearly:connected_theorem}
Let the DODS~(\ref{DODS})  
admit a four-dimensional symmetry algebra
realized by linearly connected vector fields.
Then the DODS can be transformed into  the linear DODS~(\ref{inhomogeneous}).
\end{theorem}

\noindent {\it Proof.}
Precisely two different {non}isomorphic
four-dimensional Lie algebras of linearly connected vector fields in two dimensions
exist~\cite{Gonzalez1992, DorKozMelWin2018a}.
Both were listed in Table~1 of Ref.~\cite{Gonzalez1992}
in the infinite series of algebras 20 and 21. 
Any four-dimensional Lie algebra realized by linearly connected vector fields 
can hence be transformed into one of them.

\begin{enumerate}

\item

Solvable algebra  $ {\sso {s}}_{4,3} $
with operators  $ \mbox{\bf A}_{4,5}   $ (also $ H_3$ in~(\ref{series_H}))
\begin{equation}      \label{linearly_dim4case21}
X_1 =  {\ddy}    ,
\quad
X_2 =  x  {\ddy}   ,
\quad
X_3 =   \chi   (x)   {\ddy} ,
\quad
X_4 = y  {\ddy}    ,
\quad
 \ddot{\chi}    (x)   {\not \equiv } 0 .
\end{equation}
%   where $ \{ 1 ,  x  ,   \chi   (x)    \} $  are independent functions.
In this case there are three invariants
in the space  $ ( x, y ,x_- , y_- ,\dot{y} ,   \dot{y}_-  ,  \ddot{y} )     $:
\begin{equation*}
I_1 = x,
\qquad
I_2 = x_-,
\qquad
I_3 =
{
\left|
\begin{array}{cccc}
\ddot{y}           &   \dot{y}             &    y              &   y_-              \\
0                  &   0                   &    1              &   1                \\
0                  &   1                   &    x              &   x_-              \\
\ddot{\chi} (x)  &   \dot{\chi} (x)    &    {\chi} (x)   &   {\chi} (x_-)   \\
\end{array}
\right|
\over
\left|
\begin{array}{cccc}
\dot{y}             &  \dot{y}_-             &  y                &   y_-              \\
0                   &  0                     &  1                &   1                \\
1                   &  1                     &  x                &  x_-               \\
\dot{\chi} (x)    &  \dot{\chi} (x_-)    &  {\chi} (x)     &   {\chi} (x_-)   \\
\end{array}
\right|
}    .
\end{equation*}

We obtain the invariant linear homogeneous DODE
\begin{equation}      \label{linearly_DODE_21}
\left|
\begin{array}{cccc}
\ddot{y}           &   \dot{y}             &    y              &   y_-              \\
0                  &   0                   &    1              &   1                \\
0                  &   1                   &    x              &   x_-              \\
\ddot{\chi} (x)  &   \dot{\chi} (x)    &    {\chi} (x)   &   {\chi} (x_-)   \\
\end{array}
\right|
= f (x)
\left|
\begin{array}{cccc}
\dot{y}             &  \dot{y}_-             &  y                &   y_-              \\
0                   &  0                     &  1                &   1                \\
1                   &  1                     &  x                &  x_-               \\
\dot{\chi} (x)    &  \dot{\chi} (x_-)    &  {\chi} (x)     &   {\chi} (x_-)   \\
\end{array}
\right|
\end{equation}
with the invariant delay relation
\begin{equation}      \label{linearly_delay_21}
x_-  =g (x)   .
\end{equation}
Note that the DODE is a second order one if and only if the
condition
\begin{equation*}
\left|
\begin{array}{ccc}
   0                   &    1              &   1                \\
   1                   &    x              &   x_-              \\
  \dot{\chi}  (x)    &    {\chi}  (x)   &   {\chi}   (x_-)   \\
\end{array}
\right|   _{ x_- =  g ( x )  } =  {\chi}  (x)     -  {\chi}  ( g (x)
)  -       \dot{\chi}  (x)    ( x - g ( x ) ) {\not \equiv }  0
\end{equation*}
holds.

\item

Abelian  algebra $ 4 {\sso {n}}_{1,1}  $
with operators   $ \mbox{\bf A}_{4,22}  $
(also $ S_3$ in~(\ref{series_S}))
\begin{equation}      \label{linearly_dim4case20}
X_1 =  {\ddy} ,
\qquad
X_2 = x   {\ddy}   ,
\qquad
X_3 =   \chi_2 (x)  {\ddy}  ,
\qquad
X_4 =   \chi_3 (x)   {\ddy}  ,
\end{equation}
where $ \{ 1 ,  x  ,  \chi_2 (x)   ,   \chi_3 (x)    \} $  are linearly independent functions.

In the space  $  (  x, x_-, y , y_-  , \dot{y} ,  \dot{y} _- ,  \ddot{y} ) $
there are three invariants
\begin{equation*}
I_1 = x,
\qquad
I_2 = x_-,
\qquad
I_3 =
\left|
\begin{array}{ccccc}
\ddot{y}   &   \dot{y}    &   \dot{y}_-    &   y   &   y_-    \\
0     &  0    &  0    &  1   &   1    \\
0     &  1    &  1    &  x   &   x_-  \\
\ddot{\chi}  _2  (x)  &  \dot{\chi} _2  (x)    &  \dot{\chi}  _2 (x_-)  &  {\chi} _2  (x)    &   {\chi} _2  (x_-)  \\
\ddot{\chi}  _3  (x)  &  \dot{\chi} _3  (x)    &  \dot{\chi}  _3 (x_-)  &  {\chi} _3  (x)    &   {\chi} _3  (x_-)  \\
\end{array}
\right|   .
\end{equation*}

We obtain the invariant  linear inhomogeneous DODE
\begin{equation}      \label{linearly_dim4case20_proof}
\left|
\begin{array}{ccccc}
\ddot{y}   &   \dot{y}    &   \dot{y}_-    &   y   &   y_-    \\
0     &  0    &  0    &  1   &   1    \\
0     &  1    &  1    &  x   &   x_-  \\
\ddot{\chi}  _1  (x)  &  \dot{\chi} _1  (x)    &  \dot{\chi}  _1 (x_-)  &  {\chi} _1  (x)    &   {\chi} _1  (x_-)  \\
\ddot{\chi}  _2  (x)  &  \dot{\chi} _2  (x)    &  \dot{\chi}  _2 (x_-)  &  {\chi} _2  (x)    &   {\chi} _2  (x_-)  \\
\end{array}
\right|   = f (x)
\end{equation}
with the same invariant delay relation~(\ref{linearly_delay_21}). 
% \begin{equation}      \label{linearly_delay_20}
% x_-  =g (x)   .
% \end{equation}
The obtained DODE is of order 2 if and only if
\begin{equation*}
\left|
\begin{array}{cccc}
  0    &  0    &  1   &   1    \\
  1    &  1    &  x   &   x_-  \\
\dot{\chi} _1  (x)    &  \dot{\chi}  _1 (x_-)  &  {\chi} _1  (x)    &   {\chi} _1  (x_-)  \\
\dot{\chi} _2  (x)    &  \dot{\chi}  _2 (x_-)  &  {\chi} _2  (x)    &   {\chi} _2  (x_-)  \\
\end{array}
\right|   _{ x_- =  g ( x )  } {\not \equiv }   0
\end{equation*}
holds.

\end{enumerate}

In both cases we obtain the statement of the theorem provided that
we have a second-order DODE. \hfill $\Box$

\medskip

Thus, DODSs which admit four linearly connected symmetries are linearizable.
% Both algebras are {sub}algebras of the algebra~(\ref{general_linear_operators}). 
Theorem~\ref{linearly:connected_theorem} tells us that
any Lie point symmetry algebra containing a four-dimensional sub{algebra} realized
by linearly connected vector fields will correspond to  an invariant DODS
that can be transformed into the linear DODE with a
solution-independent delay relation. The larger Lie algebra will at
most put constraints on the functions involved in the DODS.

In the rest of this section we  consider the symmetry properties
of the linear DODS~(\ref{inhomogeneous}).
First of all we mention a property that inhomogeneous  DODSs
share with ordinary inhomogeneous  ODSs (without delay).
Namely, if we know at least one solution  $ \sigma (x) $ of the inhomogeneous  equation
we can reduce the equation to a homogeneous one by the transformation
\begin{equation}
\bar{x} = x,
\qquad
\bar{y} = y - \sigma (x),
\end{equation}
For linear DODSs this does not affect the delay equation. 
We can formulate an important theorem.

\begin{theorem}       \label{general_symmetries}
The linear homogeneous  DODS~(\ref{homogeneous})  admits an infinite dimensional
Lie point symmetry algebra represented by the vector fields
\begin{equation}    \label{always_homogeneous}
X (\rho)  =  \rho  (x)  { \ddy} ,
\qquad
Y   =  y    { \ddy}  ,
\end{equation}
where $ \rho (x) $  is the general solution
of the homogeneous linear DODS~(\ref{homogeneous}).
\end{theorem}

\noindent {\it Proof.}
Application of the prolongation~(\ref{prolongation}) of the vector field~(\ref{operator1})
% \begin{equation*}
% X = \xi (x,y) { \ddx} +  \eta (x,y) { \ddy}
% \end{equation*}
to the delay equation in~(\ref{homogeneous}) gives
\begin{equation}
 \xi (x_- ,y_- ) = \dot{g} (x)  \xi (x,y)  .
\end{equation}
Since $x$, $y$, $y_-$,  $\dot{y}$ and $\dot{y}_-$   can be considered as independent
while $x_-$ and $\ddot{y}$ are related to them via equations~(\ref{homogeneous})
we get
\begin{equation}   \label{conncetion1}
\xi = \xi(x),  
\qquad  
\xi ( g ( x )  ) = \dot{g} (x)  \xi (x) .
\end{equation}
Now we apply the prolongation   to the DODE in~(\ref{homogeneous})
on the solutions of equation~(\ref{homogeneous}):
\begin{multline}   \label{deter2}
\eta  _{xx}  (x,y)
+ \left( 2 \eta _{xy} (x,y) - \ddot{\xi} (x) \right)    \dot{y}
+  \eta _{yy} (x,y)   \dot{y} ^2  \\
+  \left( \eta _y (x,y) - 2 \dot{\xi} (x) \right)
\left(    a_1  ( x  )  \dot{y}  +   a_2   ( x  )  \dot{y} _-  +  a_3   ( x  )  y +  a_4   ( x  )  y_-        \right)
\\
=  \xi (x)  \left(    \dot{a}_1  ( x  )  \dot{y}  +    \dot{a}_2   ( x  )  \dot{y} _-
+   \dot{a}_3   ( x  )  y +   \dot{a}_4   ( x  )  y_-       \right)
\\
+  a_1 (x) \left[  \eta  _{x}  (x,y)  +   \left( \eta _{y} (x,y) - \dot{\xi} (x) \right)    \dot{y}  \right]
\\
+ a_2 (x) \left[  \eta  _{x_-}  (x_-,y_-)  +  \left( \eta _{y_-} (x_-,y_-) - \dot{\xi} (x_-)  \right)    \dot{y}_-  \right]
\\
+ a_3 (x)   \eta    (x,y)
+ a_4 (x)   \eta    (x_-,y_-)       .
\end{multline}

Splitting equation~(\ref{deter2})
for  terms with $\dot{y}^2 $, $\dot{y} $, $\dot{y} _-$   and terms without derivatives,
we obtain the equations
\begin{equation}     \label{split1_a}
\eta _{yy} (x , y )  = 0    ,
\end{equation}
\begin{multline}      \label{split1_b}
 2 \eta _{xy} (x,y) - \ddot{\xi} (x)
+   a_1  ( x  )  \left( \eta _y (x,y) - 2 \dot{\xi} (x) \right)
\\
=     \dot{a}_1  ( x  )   \xi (x)  +  a_1 (x)  \left( \eta _{y} (x,y) - \dot{\xi} (x)  \right)     ,
\end{multline}
\begin{equation}      \label{split1_c}
  a_2   ( x  )     \left( \eta _y (x,y) - 2 \dot{\xi} (x)  \right)
=      \dot{a}_2   ( x  )     \xi (x)
+ a_2 (x)    \left( \eta _{y_-} (x_-,y_-) - \dot{\xi} (x_-)   \right)       ,
\end{equation}
\begin{multline}       \label{split1_d}
\eta  _{xx}  (x,y)
+  \left( \eta _y (x,y) - 2 \dot{\xi} (x)   \right)
\left(      a_3   ( x  )  y +  a_4   ( x  )  y_-       \right)
\\
=  \xi (x)  \left(     \dot{a}_3   ( x  )  y +   \dot{a}_4   ( x  )  y_-    \right)
\\
+  a_1 (x)   \eta  _{x}  (x,y)
+ a_2 (x)    \eta  _{x_-}  (x_-,y_-)
+ a_3 (x)   \eta    (x,y)
+ a_4 (x)   \eta    (x_-,y_-)       .
\end{multline}

From equation~(\ref{split1_a}) we obtain
\begin{equation}
\eta (x,y) = A (x) y + B (x)  .
\end{equation}
Substituting this result into equations~(\ref{split1_b})-(\ref{split1_d}),
we get the system
\begin{equation}      \label{split2_b}
 2 \dot{A} (x)  - \ddot{\xi} (x)
+   a_1  ( x  )  \left(  A(x)   - 2 \dot{\xi} (x)  \right)
=     \dot{a}_1  ( x  )     \xi (x)    +  a_1 (x)  \left(  A(x)  - \dot{\xi} (x)   \right)     ,
\end{equation}
\begin{equation}      \label{split2_c}
  a_2   ( x  )     \left( A (x) - 2 \dot{\xi} (x)  \right)
=       \dot{a}_2   ( x  )        \xi (x)
+ a_2 (x)    \left( A (x_- )   -  \dot{\xi} (x_-)   \right)       ,
\end{equation}
\begin{multline}       \label{split2_d}
\ddot{A} (x) y   +  \ddot{B} (x)
+  \left( A(x)  - 2 \dot{\xi} (x)  \right)
\left (      a_3   ( x  )  y +  a_4   ( x  )  y_-       \right)
\\
=  \xi (x)  \left(     \dot{a}_3   ( x  )  y +   \dot{a}_4   ( x  )  y_-     \right)
+  a_1 (x)   \left(    \dot{A} (x) y   +  \dot{B} (x)   \right)
+ a_2 (x)    \left(    \dot{A} (x_-) y_-   +  \dot{B} (x_-)   \right)
\\
+ a_3 (x)     \left(    A (x) y   +  B (x)   \right)
+ a_4 (x)    \left(    A (x_-) y_-   +  B (x_-)   \right)       .
\end{multline}
The last equation splits for terms with $y$, $y_-$  and the remaining terms as
\begin{equation}       \label{split3_d1}
\ddot{A} (x)
+   a_3   ( x  )    \left( A(x)  - 2 \dot{\xi} (x)   \right)
=  \dot{a}_3   ( x  )      \xi (x)
+  a_1 (x)    \dot{A} (x)
+ a_3 (x)      A (x)       ,
\end{equation}
\begin{equation}       \label{split3_d2}
 a_4   ( x  )   \left( A(x)  - 2 \dot{\xi} (x)  \right)
   =       \dot{a}_4   ( x  )     \xi (x)
+ a_2 (x)      \dot{A} (x_-)
+ a_4 (x)     A (x_-)           ,
\end{equation}
\begin{equation}      \label{split3_d3}
  \ddot{B} (x)
=
   a_1 (x)   \dot{B} (x)
+ a_2 (x)   \dot{B} (x_-)
+ a_3 (x)   B (x)
+ a_4 (x)   B (x_-)         .
\end{equation}

The equation~(\ref{split2_b})  can be integrated to obtain
\begin{equation}      \label{split3_b}
 A (x)
= {  1 \over 2 }   \dot{\xi} (x)  
+  {  a_1 (x)   \over 2 }  \xi  (x) 
+ A_0     ,
\qquad
A _0 =  \mbox{const}   .
\end{equation}

Substituting~(\ref{split3_b}) into the equations~(\ref{split2_c}),    (\ref{split3_d1})  and~(\ref{split3_d2}), 
we obtain three equations for $ \xi(x) $:  
\begin{equation}      \label{condition_A}
 a_2 (x)  \dot{\xi} (x)   
+ 
\left[
 \dot{a}_2  (x)   
+ 
{ a_2 (x) \over 2 }
\left(  
- { a_1 (x)  }
+   {   a_1 (g(x))     } \dot{g} (x)
-     {    \ddot{g}   (x)   \over   \dot{g}   (x)   } 
\right)
\right]
 \xi  (x)  
=  0    , 
\end{equation}
\begin{equation}      \label{condition_B}
\dddot{\xi} (x)  
+ 
\left[ 
2 \dot{a}_1 (x)   - ( a_1(x) )  ^2   - 4 a_3 (x)  
\right]
\dot{\xi} (x)  
+ 
\left[ 
\ddot{a}_1 (x)    - a_1 (x)   \dot{a} _1 (x)    - 2 \dot{a}_3 (x)   
\right]
\xi  (x)  
= 0   , 
\end{equation}
\begin{multline}    \label{condition_C}
 { a_2 (x)  \over    \dot{g}   (x)   }   \ddot{\xi} (x)  
+ 
\left[
 a_2 (x)  
\left(   
{ a_1 ( g(x)) } 
+ 
{    \ddot{g}   (x)   \over  \dot{g}   (x) ^2    } 
\right)
+ 4 a_4 ( x ) 
\right]
\dot{\xi} (x)  
\\
+ 
\left[
 2 \dot{a}_4 (x) 
+  a_2 ( x) 
\left( 
 { \dot{a}_1 (g(x))  }   \dot{g}   (x) 
+ {  a_1 (g(x)) }    {    \ddot{g}   (x)   \over   \dot{g}   (x)   } 
 + {     \dddot{g}   (x)   \over   \dot{g}   (x)  ^2  } 
-   {    \ddot{g}   (x)   ^2 \over  \dot{g}   (x)   ^3 } 
\right)
\right.
\\
\left.
+ a_4 (x) 
\left(
- {  a_1 (x)  } 
+ {  a_1 (g(x))  }    \dot{g}   (x)  
+  {    \ddot{g}   (x)   \over  \dot{g}   (x)   }  
\right) 
\right]
{\xi} (x)  
= 0    .
\end{multline}
Here in addition to Eq.~(\ref{split3_b})  we used Eq.~(\ref{conncetion1}), 
differential consequences of  both these equations 
as well as the delay equation $ x_- = g(x) $.

Three differential equations    (\ref{condition_A}),(\ref{condition_B}),(\ref{condition_C})   
together with the discrete equation~(\ref{conncetion1}) 
form an overdetermined differential-discrete  system.  
For general $  a_i (x)  $ and $g(x)$ it has only one solution, 
namely  $ \xi (x) \equiv  0 $.

Thus  we obtain symmetries with coefficients
\begin{equation}      \label{symmeties_xi_0}
\xi (x,y)   \equiv    0 ,
\qquad
\eta (x,y) =  A_0 y  + B (x)    ,
\qquad
A_0 = \mbox{const} , 
\end{equation}
 where $ B(x) $ solves 
the considered homogeneous linear DODS~(\ref{homogeneous}). 
% \begin{multline}      \label{split4_d3}
% \ddot{B} (x)
% =  a_1 (x)   \dot{B} (x)
% + a_2 (x)   \dot{B} (x_-)
% + a_3 (x)   B (x)
% + a_4 (x)   B (x_-)      ,
% 
% x_- = g (x) .
% \end{multline}
 Finally, we can rewrite the admitted symmetries
as  given in~(\ref{always_homogeneous}).
\hfill $\Box$

\medskip

\begin{theorem}   \label{inhomogeneous_symmetries} 
For the inhomogeneous linear DODS~(\ref{inhomogeneous})
the theorem is virtually the same.
The particular solution $ \sigma (x) $ 
of the inhomogeneous DODS~(\ref{inhomogeneous})   
enters into the definition of $ Y $ and we have
\begin{equation}    \label{always_inhomogeneous}
 X  (\rho)  =  \rho  (x)  { \ddy} ,
\qquad
 Y  (\sigma ) =   ( y  - \sigma (x))   { \ddy}
\end{equation}
instead of~(\ref{always_homogeneous})
\end{theorem}

The proof is obvious. Moreover, this is a property shared by linear DODSs 
and linear ODEs (without delay).

The four-dimensional Lie algebras~(\ref{linearly_dim4case21}) 
and~(\ref{linearly_dim4case20}) of Theorem~\ref{linearly:connected_theorem} 
are {sub}algebras of the infinite-dimensional Lie algebras~(\ref{always_homogeneous}) 
of Theorem~\ref{general_symmetries} 
or algebras~(\ref{always_inhomogeneous}) 
of Theorem~\ref{inhomogeneous_symmetries}. 
% The symmetry group corresponding to~(\ref{always_homogeneous})  
% of Theorem~\ref{linearly:connected_theorem}  
% is present for all linear DODSs. 
% The algebra is solvable, its {nil}radical $ X ( \rho )$ is Abelian. 
Theorems~\ref{general_symmetries} and~\ref{inhomogeneous_symmetries} 
are very useful for recognizing linearizable DODSs. 
They can however not be used for symmetry reduction so they do not provide 
any explicit analytic solutions. 
Indeed, knowing all symmetries specified in~(\ref{always_homogeneous})  
in the case of homogenous DODSs  
(in~(\ref{always_inhomogeneous})   in the case of inhomogenous DODSs)  
is equivalent to knowing all solutions.

The important question that remains is does the linear homogeneous DODS~(\ref{homogeneous}) 
allow any further symmetries for special cases of the functions $ a_i (x) $  and $ g (x) $? 
The answer is given by Theorem~\ref{more_symmetries_1} below.

Now we turn to the examination of the equations~(\ref{conncetion1}),  
(\ref{condition_A}),    (\ref{condition_B}) and~(\ref{condition_C})  
which $ \xi (x) $ must satisfy. 
We will make use of the following result.

\begin{lemma}
Let us consider an overdetermined differential-discrete system 
\begin{subequations}     \label{dd_system}
\begin{gather}  
\xi ( g ( x )  ) = \dot{g} (x)  \xi (x)  ,     \label{d1_system}  \\
\dot{\xi} (x)    = K (x)     \xi  (x)   .       \label{d2_system}  
\end{gather}  
\end{subequations}  
{Non}trivial solutions exist if and only if 
the functions $g(x)$ and $K(x)$ satisfy the compatibility condition
\begin{equation}        \label{compatibil}
K(g(x)) ( \dot{g} (x) ) ^ 2  =  \ddot{g} (x)  + K(x)  \dot{g} (x)
\end{equation}
These solutions can be presented as 
\begin{equation}    \label{Ksolution}
\xi  (x)    =  \xi  _0     e ^{ \int _0 ^x      K(s)  ds }   , 
\qquad 
  \xi  _0 \neq 0  .  
\end{equation}
\end{lemma}

\noindent {\it Proof.}
Differentiating~(\ref{d1_system})  with respect to $x$ we get 
\begin{equation*}
 \dot{\xi} ( g ( x )  ) \dot{g} (x)  = \ddot{g} (x)  \xi (x)  + \dot{g} (x)  \dot{\xi} (x) .
\end{equation*}
Using both equations of the considered system~(\ref{dd_system}),  
we obtain an identity 
\begin{equation*}
\left[   K(g(x)) ( \dot{g} (x) )^2  -  \ddot{g} (x)  -  \dot{g} (x)  K(x)   \right] \xi  (x)    =  0  .
\end{equation*}
We conclude that either $ \xi (x)  \equiv 0 $ or the  condition~(\ref{compatibil}) must hold. 
The function   $ \xi ( x ) $  must still satisfy the first order ODE~(\ref{d2_system}). 
We can write the general solution as given in~(\ref{Ksolution}). 
\hfill $\Box$

\medskip

Using this lemma, we can establish the following result.

\begin{theorem}    \label{more_symmetries_1}
Let the linear homogeneous DODS~(\ref{homogeneous})  
have functions $ a_i (x) $  and $ g (x) $ such that 
the equations~(\ref{conncetion1}),    (\ref{condition_A}),      (\ref{condition_B})    and~(\ref{condition_C}) 
have a {non}trivial solution $ \xi(x) $, 
then the symmetry algebra is larger. 
In addition to the symmetries~(\ref{always_homogeneous}) 
it contains one additional symmetry of the form  
\begin{equation}    \label{one_more}
Z 
=   \xi (x) {\ddx} 
+ {  1 \over 2 } 
\left( 
   \dot{\xi} (x)  
+  {  a_1 (x)  }  \xi  (x) 
\right) 
y {\ddy}   . 
\end{equation} 
\end{theorem}

\noindent {\it Proof.} We have to consider two cases

\begin{enumerate}

\item

$ a_2 (x) { \not \equiv } 0 $

In this case we rewrite  the equation~(\ref{condition_A}) as 
\begin{equation}      \label{split3_c_a}
\dot{\xi} (x)
= K_1 (x)    \xi  (x) , 
\qquad 
K_1 (x)  
= 
-  { \dot{a}_2  (x)    \over   a _2  (x)    }
+ { a_1 (x)  \over 2 }
-  { a_1 (g(x))  \dot{g} (x)    \over 2 }
+  {     \ddot{g}   (x)   \over 2   \dot{g}   (x)   }   .
\end{equation}
% Here  we used
% \begin{equation*}
%  \dot{\xi}    ( g ( x )  ) =   { \ddot{g} (x)   \over  \dot{g} (x)   }  \xi (x)   +   \dot{\xi} (x)  ,
% \end{equation*}
% obtained by differentiation of~(\ref{conncetion1}).

The equations~(\ref{conncetion1}) and~(\ref{split3_c_a}) 
can have a {non}trivial solution $ \xi (x)$ if the compatibility condition~(\ref{compatibil}) holds. 
If in addition this solution $\xi(x)$ and functions  $  a_i (x)  $ and $g(x)$ 
satisfy the equations~(\ref{condition_B}) and~(\ref{condition_C}), 
then there is an additional symmetry of the form~(\ref{one_more}) with 
\begin{equation}   \label{additional_symmetry_1}  
\xi (x) =    e ^{ \int       K _1 (x)  dx }   . 
\end{equation}

\item

$ a_2 (x)  \equiv   0 $,  $ a_4 (x) { \not \equiv } 0 $

Now the equation~(\ref{condition_A})  holds identically. 
The equation~(\ref{condition_C})  can be rewritten as 
\begin{equation}      \label{split3_c_b}
\dot{\xi} (x)
= K_2 (x)    \xi  (x) , 
\qquad 
K_2 (x)  
= 
- { \dot{a}_4 (x)  \over 2  a_4 (x)   }
+ { a_1 (x)  \over 4 }
- { a_1 (g(x))  \dot{g} (x)    \over 4 }  
- {     \ddot{g}   (x)   \over 4   \dot{g}   (x)   }    .
\end{equation}
The overdetermined differential--discrete system~(\ref{conncetion1}) and~(\ref{split3_c_b})
has a {non}trivial solution if functions  $g (x)$ and  $ K_2  (x) $ satisfy 
the compatibility  condition~(\ref{compatibil}). 
If in addition  this solution $ \xi(x) $ and functions $ a_1 (x) $  and  $ a_3 (x) $  
satisfy  the equation~(\ref{condition_B}), then there is an additional symmetry~(\ref{one_more}) 
with 
\begin{equation}   \label{additional_symmetry_2}  
\xi (x) =    e ^{ \int       K _2 (x)  dx }   . 
\end{equation}   
\end{enumerate}
\hfill $\Box$

\medskip

The additional symmetry~(\ref{one_more})   can be used to simplify the DODSs. 
For convenience we formulate this result only for the linear homogenous DODSs. 
However,  the same result can be established for   the linear inhomogenous DODSs.

\begin{theorem}    \label{more_symmetries_2}
Let the linear homogeneous DODS~(\ref{homogeneous})  
allow   the additional symmetry~(\ref{one_more}). 
Then the DODS can be transformed to 
\begin{equation}      \label{canonical1}
\ddot{y} =  \alpha  \dot{y}_-   + \beta   y   +  \gamma    y_-    ,
\qquad
x_- = x - C
\end{equation}
with
\begin{equation*}
\alpha ^2 + \gamma^2   \neq 0 ,
\qquad
C  >  0  , 
\end{equation*}
where   $ \alpha  $,  $ \beta $,   $ \gamma $ and $ C $  are constants. 
This DODS admits the symmetry algebra 
\begin{equation}      \label{canonical1_symmetries}
 X  ( \rho )  =  \rho  (x)  { \ddy} ,
\qquad
 Y   =   y    { \ddy}  ,
\qquad
 Z   = { \ddx}   . 
\end{equation} 
\end{theorem}

\noindent {\it Proof.} 
Let us consider a linear homogenous DODS~(\ref{homogeneous})
which admits a symmetry of the form~(\ref{one_more}) with
$ \xi (x) {\not \equiv} 0 $.
% First we note that we can always transform  the inhomogeneous DODS~(\ref{inhomogeneous})
% into its homogeneous counterpart~(\ref{homogeneous}),
% which admits symmetries~(\ref{addhomoalways}).

Let us find a change of variables 
$ (x,y) \rightarrow ( \bar{x}, \bar{y}) $  
that straightens out the vector field  $Z$ of (\ref{one_more}) 
into $ \bar{Z} = \partial _{\bar{x}}  $ 
and preserves the linearity of the DODS. 
To preserve the linearity the transformation must have the form 
\begin{equation} 
 \bar{x}  = A(x) ,  
\qquad 
 \bar{y} = B(x) y + C(x) , 
\qquad 
B(x)  \neq 0 .  
\end{equation} 
In the new variables we have 
\begin{equation} 
\bar{Z} = \xi \dot{A} { \partial  \over \partial \bar{x} }
+ 
\left( 
\left(
\xi \dot{B} 
+    
{1 \over 2 }  \left( 
  \dot{\xi} (x)  
+    a_1 (x)     \xi  (x) 
\right)   B   
\right)  y  
+ \xi \dot{C} 
\right)
{ \partial  \over \partial \bar{y} }
\end{equation} 
and we impose
\begin{equation} 
 \xi \dot{A} = 1 , 
\qquad 
\xi \dot{B} 
+    
{1 \over 2 }  \left( 
  \dot{\xi} (x)  
+    a_1 (x)     \xi  (x) 
\right)   B   
= 0 , 
\qquad 
\dot{C} = 0    . 
\end{equation} 
The solution of this system is 
\begin{equation}
A =  \int   { 1 \over \xi (x) } dx   + A_0  , 
\qquad 
B = { B_0     \over \sqrt{\xi} }  e^{ -  { 1 \over 2}  \int {a}_1 (x) dx } y  , 
\qquad 
C = C_0 , 
\end{equation}
where  $ A_0 $, $ B_0 $ and $ C_0 $ are inessential constants. 
Choosing $ A_0 = 0 $, $ B_0 =1$ and $ C_0  = 0 $, 
we obtain the new variables 
\begin{equation}
\bar{x}  =  \int   { 1 \over \xi (x) } dx   , 
\qquad 
\bar{y}  = { 1      \over \sqrt{\xi} }  e^{ -  { 1 \over 2}  \int {a}_1 (x) dx } y   . 
\end{equation}
The DODS that we obtain   is linear and invariant under transformation of $ \bar{x} $. 
Dropping the bars on  $ \bar{x} $ and $ \bar{y} $, we rewrite the DODS as 
\begin{equation}        \label{transformed} 
\ddot{y} =   a_1  \dot{y}  +  a_2   \dot{y}_-   + a_3   y   +  a_4   y_-    ,
\qquad
x_- = x - C , 
\end{equation}
where $ a_1 $ and $C$  are constants.    
The translational invariance corresponding to $ Z = \partial _x $ 
would allow $ a_i$  to be functions of $ \Delta x = x- x_- $, 
however the delay condition imposes     $ \Delta x = \mbox{const}$.

A further transformation 
\begin{equation}
 \bar{y} = e^{ -  { 1 \over 2}   \tilde{a}_1  x }   y
\end{equation}
takes~(\ref{transformed})  
into~(\ref{canonical1}) 
with symmetry algebra~(\ref{canonical1_symmetries}). 
\hfill $\Box$

\medskip

\section{Group invariant solutions of DODSs}

\label{Invariant solutions}

% \subsubsection{General method}

% The main reason for applying Lie  group and Lie algebra theory 
% for equations of any type is to obtain exact solutions. 
% For ODEs (without delay) a {non}trivial symmetry makes it possible to reduce the order of the ODE. 
% In particular, if the dimension of the symmetry group is large enough, 
% it is possible to reduce the order to zero, 
% i.e. obtain the general solution in explicit or at least implicit form.  
% {\bf 
% It is also possible to look for solution of a particular form, 
% which are invariant with respect to the Lie group action 
% generated by the considered symmetries. 
% }
% For PDEs it is  possible to decrease the number of independent variables 
% and thus obtain particular solutions satisfying specific boundary or initial conditions. 
% In both cases this "symmetry reduction" is achieved by looking for solutions 
% that are invariant under s subgroup of the symmetry group of the equation. 

% In our previous article~\cite{DorKozMelWin2018a} we showed that invariance 
% can be used to reduce a first-order DODE with one delay point 
% to two functional equations. 
% Here we will adapt the method to the case of second-order DODS (with one delay point $ x_-$). 

The main reason for applying Lie group and Lie algebra theory to equations of any type 
is to obtain exact solutions. 
For PDEs, once a symmetry group is known, 
one looks for solutions invariant under a subgroup $G_0 \subseteq G $, 
i.e. solutions annihilated by the Lie algebra    $ L_0$ of $G_0$. 
This amounts to adding one, or more linear first-order PDEs to the one which   we are solving. 
By first solving  these additional equations using the method of characteristics 
we obtain the general form of the solution $y ( {\bf x} )$,   $ {\bf x} \in \mathbb{R}  ^p$ 
in terms of the invariants  $ I _k $, $ k< p$ of the subgroup $G_0$.  
The reduced solution depends on fewer variables since we have $ k< p$. 
Solving the reduced equation, we obtain a subset of solutions satisfying specific boundary or initial 
conditions (see e.g. in Olver's book~\cite{Olver1986}).

For ODEs it is also possible to obtain group invariant solutions 
by imposing invariance under a subgroup of the symmetry group. 
This again provides very specific solutions for a restricted type of initial conditions. 
A more efficient procedure for ODEs is to take  
an element of the symmetry algebra and 
to "straighten out" the corresponding vector field 
by a change of variables
\begin{equation} 
(x,y) \rightarrow  (t,u) , 
\qquad 
X = \xi (x,y) {\ddx}  +  \eta  (x,y) {\ddy} 
\rightarrow
\tilde{X} = {\partial  \over \partial u} . 
\end{equation} 
The same transformation will take the original ODE into one not containing the variable $u$
 \begin{equation} 
E ( x,y,  \dot{y},  \ddot{y}, ... ,  y ^{(n)} ) = 0   
\rightarrow
\tilde{E} ( t , \dot{u},  \ddot{u}, ... , u ^{(n)} ) = 0    . 
\end{equation} 
Putting $   \dot{u} = v (t) $, we obtain an ODE for $v(t)$  
of order $ n-1 $. 
If we can integrate it to obtain $ v = v ( t, C_1, ..., C_{n-1})$, 
then we can integrate $v$ to obtain $ u = u ( t,   C_1, ..., C_{n})$, 
i.e. the general solution. 
If $\tilde{E}= 0 $ allows a symmetry group, we can repeat the procedure and 
reduce the order to $ n-2 $.

From this point of view DODSs are more similar to PDEs than to ODEs.     
In our earlier article~\cite{DorKozMelWin2018a} 
we used a one-dimensional {sub}algebra 
$ L_0     \subseteq  L  $ of the symmetry algebra of the given DODS 
to obtain group invariant solutions. 
These were particular solutions satisfying very specific initial conditions.  
The same method works for DODSs of any order (with one delay).

Thus, we are given a second-order DODS~(\ref{DODS}), 
i.e. the functions $f$  and $g$ are given. 
We calculate the symmetry algebra of the DODS, i.e. solve the equations
\begin{subequations}     \label{invariance}
\begin{gather}
{\bf pr} \ X    \left[    \ddot{y} - f ( x, y, y_- , \dot{y}, \dot{y}_-  )  \right]_{ \ddot{y} = f , \ x_- = g } =  0  , 
\\
{\bf pr} \ X    \left[    x_-  - g ( x, y, y_- , \dot{y}, \dot{y}_-  )  \right]_{ \ddot{y} = f , \ x_- = g } =  0
\end{gather}
\end{subequations}
for the coefficients $ \xi (x,y)  $ and $ \eta (x,y)  $ in~(\ref{operator1}). 
Thus, the symmetry algebra is known and we can classify its one-dimensional {sub}algebras 
into conjugacy classes under the action of the group $G$ leaving the DODS~(\ref{DODS}) invariant. 
This group is also known: 
$  G \sim G_0   \triangleright  \mbox{exp} \ L $, 
where $ G_0$  is a discrete group leaving the DODS invariant 
($ \mbox{exp} \  L $ is the connected component of $G$).  
For methods of classifying {sub}algebras of Lie algebras see~\cite{Patera1977a, Patera1977b} 
and references therein.

Each {sub}algebra is represented by one vector field of the form~(\ref{operator1}). 
We must find the invariants of the corresponding one-parameter subgroup 
acting in the $ (x,y) $ plane.   
All invariants will be functions of three basic invariants 
$ I _k ( x, y, x_-, y_- ) $, $ k = 1, 2, 3 $  
determined from the equation 
\begin{equation} 
\left[
\xi (x,y) {\ddx} + \eta(x,y) {\ddy}  
+ \xi (x_-,y_-) {\partial \over  \partial x_- } + \eta(x_-,y_-) {\partial \over  \partial y_- } 
\right]
I ( x, y, x_-, y_- )   = 0   . 
\end{equation} 
For the method to provide explicit solutions 
there must exist a set of basic invariants satisfying 
\begin{equation}   \label{basic}
J_1= J_1 (x,y) ,
\qquad
J_2 = J_2  ( x, y, x_- , y_- )  , 
\qquad
J_3 = J_3 (x_-,y_-)   =  J_1 (x_-,y_-) 
\end{equation}
together with the condition 
\begin{equation}    \label{Jacobian} 
\mbox{det}   \left(  { \partial ( J_1 , J_2 )  \over  \partial ( y  , x _- ) }  \right)     \neq   0  .
\end{equation}
This allows us to set the invariants equal to constants 
\begin{equation}  
J_1 (x,y) = A,
\qquad
J_2  ( x, y, x_- , y_- )  = B ,  
\qquad
J_1 (x_-,y_-) = A 
\end{equation}
and to solve for $y$ and $x_-$  
obtaining the "reduction formulas"
\begin{equation}     \label{reduction_formulas}
y = h (x, A) ,
\qquad 
y_-  = h (x_- , A) , 
\qquad
x_-  = k ( x, A, B)  . 
\end{equation}
According to the implicit function theorem 
the functions     $ h (x, A) $ and  $ k ( x, A, B) $  
are determined uniquely up to the indicated constants, 
at least in some neighbourhood of $ x \in I$. 
Substituting  
into the DODS,  we obtain equations for the constants $ A $  and $ B $ 
and~(\ref{reduction_formulas}) 
then gives the exact invariant solution.  
Unfortunately, this method for DODSs give only particular solutions. 
Thus, the initial data on the first step   $ x \in [ x_{-1} , x_0 ] $  
must satisfy $ y = h (x, A )  $.

Examples will be given in Section~\ref{micro_model}.

Let us mention that another approach to finding exact solutions 
of delay differential equations, which makes use of functional separation of variables,  
was proposed 
in  \cite{bk:PolyaninZhurov[2014]a, bk:PolyaninZhurov[2014]b, PolyaninSorikin, PolyaninZhurov2015}

\section{Exact analytical solutions of a DODS describing traffic flow}

\label{micro_model}

Traffic flow is one of the areas in which time delay plays an important role.
The time delay $ \tau $ corresponds to the time
it takes a driver of a car to react  to changes of velocity  of a preceding vehicle in a row of cars.
One of the standard nonlinear  "follow-the-leader" models  is given
by the system of delay ordinary differential equations~\cite{Gazis1961, Nagel2003}
\begin{equation}   \label{acc2}
\ddot{x} _{n+1} (t)
=  \alpha  \left[  \dot{x} _{n+1} (t)  \right] ^{n_{1}}
{ \dot{x}_{n} (t-\tau)  - \dot{x} _{n+1} (t-\tau)   \over  \left[  x_{n}(t-\tau)   -   x _{n+1} (t-\tau) \right] ^{n_{2}} }  ,
\qquad
 \alpha   \neq 0 ,
\end{equation}
where $ \alpha  $,  $ n_1 $, $ n_2 $ are adjustable constant parameters
and $ \tau $ is the constant time delay  $ t - t_- = \tau > 0 $.
The function $ x_n (t) $ is the position of the $n$-th vehicle at time $t$,
the dots denote time derivatives.
Equations~(\ref{acc2}) can be viewed as a finite ($ 0 \leq n \leq N < \infty$)
or infinite  ($ 0 \leq n < \infty$) set of DODEs,
or as a delay ordinary differential-difference equation
(since     $ x_n (t) $ is a function of a continuous variable $t$ and a discrete one $ n $).
The leading car's position is $ x _0 (t) $
and must be given as part of the input in the model.

In order to apply the Lie theory formalism developed in this article
we restrict to the case of two cars: the leading one with position $
x _0 (t) $  and a following one with position $  x_1 (t) \equiv x
(t) $. The system~(\ref{acc2}) then simplifies to the DODE
    \begin{equation}   \label{acc}
\ddot{x}
=  \alpha    \dot{x}    ^{n_{1}}
{ \dot{x} _{0, - }     - \dot{x}  _-    \over  (   x _{0, - }    -   x   _-   ) ^{n_{2}} }  ,
\qquad
 t_- = g ( t , x , x_- , \dot{x} , \dot{x} _ -  )    ,
\end{equation}
where we replaced $ (x,y) $  by $ (t,x) $ (as usual in traffic flow
studies) and allowed $ t _-  $ to be a function of the indicated
variables, instead of imposing $ t_- = t - \tau $ with constant $
\tau $. We have put $ x  _-  = x ( t_- ) $,  $ x  _{0,-}  = x_0  (
t_- ) $. 
The function $ x  _{0} (t) $ must be given for all $t$, $ t_{-1} = t_0
- \tau < t < \infty$.

 A group analysis of equations~(\ref{acc}) and also~(\ref{acc2})
is in progress.
Here  we will just present several examples showing how group analysis can lead to exact solutions
invariant  under some subgroup of the symmetry group.

% \subsection{Example 0  (from Pavel)}

\begin{example}

Let us start with the simplest and the most studied case.
The leading car is moving with a constant velocity   $ v $, $ x_0 (t) = v t $,
and the time delay is constant.
We then have
\begin{equation}   \label{example_0}
\ddot{x}
=  \alpha    \dot{x} ^{ n_1}
{  v   -  \dot{x} _-    \over  (   v t_-    - x_- ) ^{n_2}    }  ,
\qquad
 t_-    =  t  - \tau   .
\end{equation}

This DODS is invariant under transformations corresponding to
\begin{equation}
X    =        {\partial \over \partial t }   +    v     {\partial \over \partial x }    .
\end{equation}
Invariants of $ X $ in the space $ ( t, t_-, x, x_- ) $ are
\begin{equation*}
I_1 =  x - v t   ,
\qquad
I_2 = t - t_ - ,
\qquad
I_3 = x_-   - v t_-     .
\end{equation*}
We substitute
\begin{equation}
x =  v t + A ,
\qquad
t_- = t - B
\end{equation}
into~(\ref{example_0}) and obtain  $ B =\tau  $.
Thus, we have a particular solution
\begin{equation}
x(t)  =  v t + A ,
\qquad
A < 0 ,
\qquad
t_- = t - \tau .
\end{equation}

While the above considerations are very simple and obvious,
they do indicate that the model~(\ref{acc2})   makes sense.
If the leading car moves with constant velocity,
the following car should move with the same velocity.
If we have $ A = 0 $ the two cars would be in the same position
(i.e. on top of each other if we have a single lane).
Hence we impose $ A \neq 0 $, or    $ A < 0 $
since $ x(t) $ is the position of the second car.
\end{example}

% \subsection{Examples  2--5 (not rewritten)}

\begin{example}

% \subsubsection{Example 2}

Consider the  equation~(\ref{acc}) with
\begin{equation*}
n_{1}(n_{1}-1)\neq0, \qquad n_{2}=0, \qquad x_0 (t)=k t^{n}, \qquad
n \neq 0 ,
\end{equation*}
where $n=1-n_{1}^{-1}$.
The DODS
\begin{equation}   \label{example_2}
\ddot{x}
=  \alpha     \dot{x}  ^ {n_1}
( k n   t _- ^{n-1}   -  \dot{x} _-  )    ,
\qquad
 t_-        = q  t  ,
\qquad
0 < q  <  1
\end{equation}
admits a two-dimensional symmetry algebra with basis
\begin{equation}
X_1 =  t    {\partial \over \partial t }      +  n   x {\partial
\over
\partial x }   ,
\qquad X_2 = {\partial \over \partial x } .
\end{equation}
A solution invariant under the {sub}group generated by $ X = X_1 - n \beta X_2 $ has the
form
\begin{equation}     \label{solution_form_2}
x(t)=    { \beta } +A t^{n},
\qquad
t_{-}=B t.
\end{equation}
$A$ and $B$ are integration constants, $\beta$ is the result of an
$x$ translation.

Substituting~(\ref{solution_form_2}) into the DODS
(\ref{example_2}), we get restrictions  for the constants
\begin{equation}     \label{constraint_2}
 \alpha   {  (  n_1 -  1   )  ^{n_1}   \over n_1  ^{n_1 -1} }   q ^{ - { 1 \over n_1 }  }    ( k -  A )  A ^{ n_1 - 1 }  =  - 1  ,
\qquad
B = q   .
\end{equation}

Again the model and the solution~(\ref{solution_form_2}) make sense.
This time the leading car moves with acceleration $ \ddot{x}_0  = k n (n-1) t^{n-2}$.
The second car tries to adapt with solution~(\ref{solution_form_2}).
The distance would be constant for $A=k$.
This is however not possible since $A=k$ contradicts the constraint~(\ref{constraint_2}).
To avoid a collision we must have $A<k$. For $ \alpha >0 $, $n_1 > 1$ a collision is inevitable.
\end{example}

% \subsection{Example 4}

\begin{example}

% \subsubsection{Example 4}

The  equation~(\ref{acc}) with
\begin{equation*}
n_{2}=n_{1} = n  \neq0,
\qquad
x_0 (t)=k e^{\varepsilon t}  
\end{equation*}
leads to the DODS
\begin{equation}   \label{example_4}
\ddot{x}
=  \alpha     \dot{x} ^ {n}
{   k \varepsilon e^{\varepsilon t_- }     -  \dot{x} _-   \over (  k e^{\varepsilon t_- }    - {x} _-   ) ^{n}   }   ,
\qquad
 t_-        =  t  -   \tau    ,
\end{equation}
which admits the symmetry
\begin{equation}
X =    {\partial \over \partial t } +   \varepsilon   x  {\partial \over \partial x }  .
\end{equation}
A solution invariant with respect to $X$ has the form
\begin{equation}     \label{solution_form_3}
x(t)=A e^{\varepsilon t},
\qquad
t_{-} =t - B.
\end{equation}
Substituting it into the DODS, we get a restriction for the
constants. They should satisfy the equation
\begin{equation}    \label{constraint_3}
 \alpha   \left[    {  \varepsilon  e ^{  \varepsilon \tau }     A  \over  k - A  }    \right] ^{n-1}   = 1 ,
\qquad
B = \tau   .
\end{equation}

In this case the acceleration of the first car is exponential (by assumption).
The second car adapts with the same type of exponential acceleration.
A constant distance between the 2 cars requires the identity $A=k$
but this is ruled out by the constraint~(\ref{constraint_3}).
To avoid a collision we must have $A<k$. This is only possible for $ \alpha >0$.
\end{example}

The last two examples show the usefulness of simple exact solutions
of DODS. They help identify useful limits on the parameters of the
model like $\alpha$, $n_1$ and $n_2$ in~(\ref{acc}). For other
discussions of traffic models and their solutions (with no use of
group theory) see  \cite{Whitham1990, Lighthill1955, Tyagi2009,
Newell1961, Chandler1958}.

\section{Conclusion}

\label{Conclusion}

Let us sum up the main results of this article devoted to second-order DODSs 
with one delay point $ x_-$ as presented in equations~(\ref{DODS}). 
We shall also put the results into a more general context.

\begin{enumerate} 

\item 

A second-order genuinely {non}linear  DODS of the type~(\ref{DODS}) can have 
a Lie point symmetry algebra $L$ of dimension $ \mbox{dim} \ L \leq 6$.  
Representatives of all such DODSs and their symmetry algebras 
are listed in Tables~$1$--$4$. 
If we have   $ \mbox{dim} \ L > 6 $, then the DODS is linearisable  
by an invertible transformation.

\item    

All linear or linearisable  DODSs of this type have infinite-dimensional 
Lie point symmetry algebras. 
These algebras are solvable and their {nil}radicals  are Abelian and of codimension 1 or 2. 
The {nil}radical   $ NR (L) $ is infinite-dimensional because it represents the linear  
superposition principle. 
The fact that the initial data depend on a function,   
$ y(x) = \varphi (x) $, $ x \in [ x_- , x_0] $, where $ \varphi (x)$  is an arbitrary function, 
implies that an infinite set of solutions     must exist. 
The factor algebra $ L / NR (L) $ always contains one    {non}nilpotent element 
that reduces to $ X = y \partial _y $ once the DODS is transformed 
into a linear homogeneous DODS.

\item

In special cases,   
when the overdetermined differential-discrete system of 
equations~(\ref{conncetion1}),    (\ref{condition_A}),      (\ref{condition_B})    
and~(\ref{condition_C})  has {non}trivial solutions $\xi (x) {\not \equiv} 0$, 
the factor algebra  $ L / NR (L) $   contains 
 a second {non}nilpotent element  
$ Z =\xi(x) \partial _x + \eta (x,y ) \partial _y $, $ \xi  (x) {\not \equiv 0} $. 
A further transformation of variables 
takes the DODS into a linear homogeneous DODS  with constant coefficients 
\begin{equation}    \label{conclusion_constant}
\ddot{y} = \alpha  \dot{y}_-   + \beta y + \gamma y_- , 
\qquad 
x_- = x - C  
\end{equation} 
and its symmetry algebra into 
\begin{equation} 
X ( \rho ) = \rho (x) {\ddy} , 
\qquad 
Y =  y {\ddy} , 
\qquad 
Z = {\ddx}  . 
\end{equation}

Let us apply the standard method of solving linear ODEs to the DODSs. 
We put 
\begin{equation}       
y = A e ^{ \lambda x } , 
\qquad 
y _- = A e ^{ \lambda (x-C)  } 
\end{equation} 
into~(\ref{conclusion_constant}) 
and obtain a "modified characteristic equation" for $ \lambda $
\begin{equation}        \label{lambda_equation} 
\lambda ^2 -  \alpha  \lambda   e^{ - \lambda C  }   -   \beta  -  \gamma e^{ - \lambda C }  = 0 . 
\end{equation}

Any solution $\lambda$ of the transcendental   equation~(\ref{lambda_equation}) 
gives a solution of~(\ref{conclusion_constant}). 
It is however difficult to obtain possible constraints on the parameters   
 $ \alpha  $, $ \beta $, $ \gamma $ and $  C $ 
for which an infinite set of solutions~(\ref{lambda_equation})  exists 
and the class of initial data for which this set of solutions is complete.   
For further information on delay ordinary differential equations, 
in particular linear once  with constant coefficients, 
we refer to Driver's book~\cite{Driver1977}.

\item 

In order to show that a DODS~(\ref{DODS}) 
is linearisable it is sufficient to show that it allows a four-dimensional  
symmetry algebra realized by linearly connected vector fields 
(see Theorem~\ref{linearly:connected_theorem}).

\item 

A one-dimensional {sub}algebra of the symmetry algebra of a DODS 
allows us to obtain particular invariant solutions of the DODS 
if the vector fields~(\ref{operator1}) 
representing this algebra satisfies the condition
$$ 
\xi (x,y) {\not \equiv } 0 . 
$$ 
Let the one parameter group $ \mbox{exp} ( \epsilon X ) $ 
acting on the space have  three invariants 
$ J_1 = J_1 (x,y) $,  $ J_2 = J_ 2 ( x , y  ,  x_- ,  y_- ) $,   $ J_3 = J_1 ( x_-. y_- ) $ 
satisfying (locally) condition~(\ref{Jacobian}). 
Than it is possible to express $y$, $y_-$  and $x_-$  as  functions of  $x$ 
and some constants as in~(\ref{reduction_formulas}). 
Ultimately we can also express  $ dy_- / dx_- $,  $ d^2y_- / dx_- ^2$ 
and all higher derivatives  $ d^k y_- / dx_- ^k $ in terms of $x$.  
Substituting   $ y = h ( x, A) $, $  x_- = k ( x, A, B) $,   $ y_- = h ( x_-, A) $ 
as well as all derivatives  $\dot{y} (x)$, $\dot{y} (x_-)$ and $\ddot{y} (x)$ 
into the DODS, 
we obtain  functional equations from which we determine the constants $A$  and $B$ 
compatible with the DODS.

\end{enumerate}

Some of the above results can be extended to DODSs of arbitrary order $n$.

\medskip

\noindent  \underline{Results 1.}   
An $n$-th order {non}linearisable DODS with one delay point $x_-$  
will have a finite-dimensional symmetry algebra. 
The reason is that all individual algebras that can be realized by vector fields 
in two dimensions (algebras 1-19 in Table~1 of~\cite{Gonzalez1992}) 
have dimensions $ \mbox{dim} \  L  \leq 8 $.   
Those in the series 20-28 in the same table 
have arbitrary large dimensions 
 $ \dim \  L = r + k $ with  $ 0 \leq  k \leq 3 $ and $ r  = \mbox{dim} \ L _0 $, 
where $ L_0 $ is  a maximal Abelian subalgebra of $ L $ and 
is realized by linearly connected vector fields.  
The actual value of $r$ is related to the order of the DODS 
(we have $ r = 1 $  for $ n=1 $~\cite{DorKozMelWin2018a, DorKozMelWin2018b} 
and $r=3$  for $n=2$ as shown above).   
The maximal dimensions of the symmetry algebra  are   
$ \mbox{dim} \ L  =3  $ and   $ \mbox{dim} \ L  =6  $, respectively.

\noindent  \underline{Results 2.}   
All linear and linearisable DODSs of any order will have infinite-dimensional 
symmetry algebras with the same structure as those for $n=1$ and $n=2$ 
and for the same reason.

\noindent  \underline{Results 3.}   
In order to show that a DODS of order $n$  is linearisable 
it is sufficient to show that it has $k$-dimensional Abelian {sub}algebra with $ k < \infty$. 
We conjecture that $k$ satisfies $ k = 2n$.

It is interesting to compare results concerning Lie point symmetries known for DODSs with 
those for ODEs without delay.  
First-order ODEs have infinite-dimensional symmetry algebras 
because of the existence of an integrating multiplier.     
For higher order ODEs the symmetry algebras are  
all finite-dimensional with $   \mbox{dim} \ L  = 8  $ for $ n=2$ and 
  $ \mbox{dim} \ L  = n  + 4   $ for $ n \geq 3 $. 
A solvable {sub}algebra  $L_0 \subseteq L$ 
of dimension $n_s$ can be used to reduce the order of the equation to $ n- n_s$. 
For $ n_s \geq n $ this means that we can reduce the order to zero and thus obtain the general solution 
(at least in implicit form). All of these properties are well known and can be found 
e.g. in Olver's book~\cite{Olver1986}.

For DODSs the situation is different.  
The concept of integrating multiplier has not been introduced. 
The symmetry algebra is infinite-dimensional if and only if the DODS 
(of any order) is linear or linearisable. 
This is similar to the case of linearisable 
partial differential equations~\cite{Kumei1982, Bluman1989}.

In both cases this algebra is useful for establishing linearizability, 
finding the "target"  linear equations and finding the linearizing transformation.

A further difference between DODSs and ODEs is that for DODSs 
invariant solutions   are very particular and satisfy very special initial conditions.

\section*{Acknowledgments}

 The research of VD was partly supported by research grant
No. 18-01-00890 of the Russian Foundation for Basic Research.
The research of PW was partially supported by a research grant from NSERC of Canada.

\eject

\section*{Appendix}

In Tables 1-4 we provide a classification of {non}linear invariant
DODSs by dimension of the Lie algebras of the admitted symmetries.
In each table the first three columns give 
Lie algebras and their
realizations. Column~1 provides the Lie algebra isomorphism class using the
notations of~\cite{SWbook}. For example, $ {\sso {n}}_{i,k} $ denotes the
$k$-th nilpotent Lie algebra of dimension $i$. The only
nilpotent algebras in the table are $ {\sso {n}}_{1,1} $, $ {\sso
{n}}_{3,1} $ and $ {\sso {n}}_{4,1} $. Similarly, $ {\sso {s}}_{i,k}
$ denotes the $k$-th solvable (but not nilpotent) Lie algebra of dimension $i$.
Simple Lie algebras are presented by their standard notations ($ {\sso
{sl}} ( 2, \mathbb{R}  ) $, $ {\sso {o}} ( 3,  \mathbb{R} ) $). 
Column~2 gives notations $\mbox{\bf A}_{i,k} $ for algebras in the list of
sub{algebras} of $ \mbox{diff} ( 2, \mathbb{R}  ) $. 
These notations were introduced in~\cite{DorKozMelWin2018a} 
for algebras of dimension $1$--$4$. 
As before $i$ is the dimension of the algebra. 
The numbers in brackets give the algebra notations used in Table~1 of Ref.~\cite{Gonzalez1992}. 
Column~3 contains vector fields spanning each representative algebra.

The two last columns give invariant DODEs and  invariant delay
relations.
For dimensions 1, 2 and 3 the tables contain all Lie algebras which can be realized  
as {sub}algebras of $ \mbox{diff} ( 2, \mathbb{R}  ) $. 
For dimensions 4, 5 and 6 only algebras which provide
invariant DODSs are given.

\eject

\subsection*{Table~1. Classification of {non}linear invariant DODSs.
Symmetry algebra dimensions 1 and 2.}

% \subsection*{Dimensions 1 and 2}

% \begin{sideways}
% \begin{equation*}
\begin{sideways}
$
\begin{array}{|c|c|l|c|c|}
\hline
    &   &   &   &  \\
\mbox{Lie algebra} & \mbox{Case} & \mbox{Operators} & \mbox{DODE} & \mbox{Delay relation}  \\
    &   &   &   &  \\
\hline
   &   &   &   &  \\
{\sso {n}} _{1,1}
&
\mbox{\bf A}_{1,1}  (9)
&
{ \displaystyle   X_1 = { \dy}   }
&
{ \displaystyle
\ddot{y} = f ( x , \Delta y   , \dot{y} , \dot{y} _-     )
}
&
{ \displaystyle
 x _-   =   g ( x  , \Delta y   , \dot{y} , \dot{y} _-     )
} \\
   &   &   &   &  \\
\hline
\hline
   &   &   &   &  \\
{\sso {s}} _{2,1}
&
\mbox{\bf A}_{2,1} (10)
&
{ \displaystyle
X_1 = { \dy}   ;
\
X_2 =  y  { \dy} }
&
{ \displaystyle
\ddot{y}   =      \dot{y}
f \left(  x    ,  {   \dot{y}  \over  \Delta y    }  ,  {   \dot{y} _-  \over  \Delta y    }     \right)
}
&
{ \displaystyle
 x_-    =      g \left(  x     ,   {   \dot{y}  \over  \Delta y    }  ,  {   \dot{y} _-  \over  \Delta y    }    \right)
}    \\
   &   &   &   &  \\
\cline{2-5}
   &   &   &   &  \\
&
\mbox{\bf A}_{2,2} (22)
&
{ \displaystyle
X_1 = { \dy}   ;
\
X_2 = x  {\dx}  +   y  { \dy}  }
&
 { \displaystyle
\ddot{y}   =  { 1 \over x  }     f  \left(  {y_x}   , \dot{y}   ,  \dot{y} _-   \right)
}
&
{ \displaystyle
   x_-    =    x   g  \left( {y_x}   ,   \dot{y}   ,  \dot{y} _-      \right)
}   \\
   &   &   &   &  \\
\hline
   &   &   &   &  \\
2 {\sso {n}} _{1,1}
&
\mbox{\bf A}_{2,3} (20)
&
{ \displaystyle
\left\{  X_1 = { \dy}  \right\} ,
\
\left\{  X_2 = x  { \dy}  \right\}  }
&
{ \displaystyle
\ddot{y}  =     f \left(  x   ,  \dot{y}   -   {y_x} , \dot{y}   -   \dot{y} _-       \right)
}
&
{ \displaystyle
   x_-    =     g \left(  x   , \dot{y}   -   {y_x} ,   \dot{y}   -   \dot{y} _-        \right)
}  \\
   &   &   &   &  \\
\cline{2-5}
   &   &   &   &  \\
&
\mbox{\bf A}_{2,4} (22)
&
{ \displaystyle
\left\{  X_1 = { \dx} \right\} ,
\
\left\{  X_2 = { \dy} \right\}  }
&
{ \displaystyle
 \ddot{y}  =  f  (   \Delta y ,   \dot{y}   ,  \dot{y}  _-   )
}
&
{ \displaystyle
 \Delta   x = g (   \Delta y ,  \dot{y}   ,  \dot{y}  _-    )
}  \\
   &   &   &   &  \\
\hline
\end{array}
$
\end{sideways}
% \end{equation*}
% \end{sideways}

\eject

% $   $

% \vspace{-20mm}

\subsection*{Table~2. Classification of nonlinear invariant DODSs. Symmetry algebra dimension
3.}

\begin{sideways}
$
\begin{array}{|c|c|l|c|c|}
\hline
  &  &  &  &  \\
\mbox{Lie algebra} & \mbox{Case} & \mbox{Operators}  & \mbox{DODE} & \mbox{Delay relation}  \\
  &  &  &  &  \\
\hline
  &  &  &  &  \\
{\sso {n}} _{3,1}
&
\mbox{\bf A}_{3,1} (22)
&
 { \displaystyle
X_1 = { \dy} ;
\
X_2 =  x  { \dy} ,
\
X_3 = { \dx} }
&
{ \displaystyle
 \ddot{y}  =   f \left(   \dot{y}   -   {y_x}   ,   \dot{y}  -  \dot{y}_-     \right)
}
&
{ \displaystyle
\Delta x  =  g \left(  \dot{y}   -   {y_x}   ,   \dot{y}   -  \dot{y}_-      \right)   .
}
\\
  &  &  &  &  \\
\hline
  &  &  &  &  \\
{\sso {s}} _{3,1}
&
\mbox{\bf A}_{3,2} ^a  (12)
&
\begin{array}{l}
{ \displaystyle X_1 =  { \dx} ,
\
X_2 =  { \dy}  ;   }
\\
{ \displaystyle
X_3 =   x  { \dx}  +  a  y  { \dy}  ,
\
0 < |a| \leq 1
  }
\\
\end{array}
&
{ \displaystyle
 \ddot{y}  =    | \Delta x | ^{a-2}
f
\left(
{     \dot{y}    \over     | \Delta x | ^{a-1}      }    ,
{    \dot{y}  _-  \over    | \Delta x | ^{a-1}     }
\right)
}
&
{ \displaystyle
 \Delta  x      =   | \Delta y | ^{1 \over  a}
g
\left(
{     \dot{y}    \over     | \Delta x | ^{a-1}      }    ,
{    \dot{y}  _-  \over    | \Delta x | ^{a-1}     }
\right)
}
\\
  &  &  &  &  \\
\cline{2-5}
  &  &  &  &  \\
&
\mbox{\bf A}_{3,3} ^a  (21, 22)
&
\begin{array}{l}
{ \displaystyle
X_1 = { \dy} ,
\
X_2 = x { \dy} ;  }
\\
{  \displaystyle
X_3 =  (1-a) x { \dx}  +  y { \dy} , }
\\
{  \displaystyle
0 < |a| \leq 1
  }
\\
\end{array}
\qquad
\mbox{i)} \   a \neq  1
&
% a \neq  1
% \qquad
 \ddot{y}  =   |x|^{  2a - 1  \over 1 -  a  }
f
\left(
{    { \displaystyle    \dot{y}   -   {y_x}   }  \over  { \displaystyle     |x|^{a \over 1 - a }   }   }  ,
{    { \displaystyle    \dot{y}  -   \dot{y} _-   }  \over { \displaystyle   |x|^{a \over 1 - a }  }  }
\right)
&
 x_-     =  x
g
\left(
{    { \displaystyle    \dot{y}   -   {y_x}   }  \over  { \displaystyle     |x|^{a \over 1 - a }   }   }  ,
{    { \displaystyle    \dot{y}  -   \dot{y} _-  }   \over   { \displaystyle    |x|^{a \over 1 - a }   }  }
\right)
\\
  &  &  &  &  \\
\cline{4-5}
  &  &  &  &  \\
&
&
\qquad
\qquad
\qquad
\qquad
\qquad
\qquad
\quad
\mbox{ii)} \
a = 1
&
% a = 1
% \qquad
 \ddot{y}  =  \left(  {    \dot{y}   -   {y_x}  }  \right)
f
\left(
x,
{   { \displaystyle   \dot{y} -  \dot{y} _-     }
\over
{ \displaystyle  \dot{y}   -   {y_x}  }  }
\right)
&
 x_-   = g
\left(
x,
{  { \displaystyle   \dot{y}  -  \dot{y} _-    }
\over
{ \displaystyle  \dot{y}   -   {y_x}  }  }
\right)
\\
  &  &  &  &  \\
\hline
  &  &  &  &  \\
{\sso {s}} _{3,2}
&
\mbox{\bf A}_{3,4}   (25)
&
\begin{array}{l}
{ \displaystyle
X_1 = { \dx} ,
\
X_2 =  { \dy} ;
  }
\\
{ \displaystyle
X_3 =   x  { \dx}  +   ( x  +  y )  { \dy}   }
\\
\end{array}
&
 \ddot{y}  =  { \displaystyle   { 1  \over   \Delta x  }  }
f
\left(
\dot{y}   -  {y_x}  ,
 \dot{y}   -   \dot{y} _-
\right)
&
 \Delta  x    =  e^ {y_x}
g
\left(
\dot{y}   -  {y_x}  ,
 \dot{y}   -   \dot{y} _-
\right)
\\
  &  &  &  &  \\
\cline{2-5}
  &  &  &  &  \\
&
\mbox{\bf A}_{3,5}  (22)
&
 { \displaystyle
X_1 = { \dy} ,
\
X_2 = x { \dy}   ;
\
X_3 =  { \dx}  +  y { \dy}
}
&
  \ddot{y} =  e^{x}   f
\left(
e^{ - x}   \left(    \dot{y}   -   {y_x}   \right)    ,
 e^{ - x}   \left(    \dot{y}   -    \dot{y} _-      \right)
\right)
&
\Delta x  =  g
\left(
e^{ - x}   \left(    \dot{y}   -   {y_x}   \right)   ,
 e^{ - x}   \left(    \dot{y}   -    \dot{y} _-      \right)
\right)
\\
  &  &  &  &  \\
\hline
\end{array}
$
\end{sideways}

\eject

\begin{sideways}
$
\begin{array}{|c|c|l|c|c|}
\hline
  &  &  &  &  \\
\mbox{Lie algebra} & \mbox{Case} & \mbox{Operators}  & \mbox{DODE} & \mbox{Delay relation}  \\
  &  &  &  &  \\
\hline
  &  &  &  &  \\
{\sso {s}} _{3,3}
&
\mbox{\bf A}_{3,6}  ^b  (1)
&
\begin{array}{l}
 {\displaystyle
X_1 =  { \dx} ,
\
X_2 =  { \dy} ; }
\\
 {\displaystyle
X_3 =   ( b x + y )  { \dx}  +  ( b  y - x )   { \dy}  ,  }
\\
 {\displaystyle
b \geq 0 } \\
\end{array}
&
\begin{array}{c}
{\displaystyle
   \ddot{y}
=  ( 1  +   \dot{y} ^2 ) ^{3/2} e^{ b   \arctan  \dot{y}   }   }
\\
\\
\times
{ \displaystyle
f
\left(
 { \dot{y}    -  {y_x}
\over
1 +    \dot{y}   {y_x}     }  ,
 {  \dot{y}   -  \dot{y}  _-
\over
1 +    \dot{y}   \dot{y}  _-        }
\right)  }
\\
\end{array}
&
\begin{array}{c}
{\displaystyle
 \Delta x     e^{ b   \arctan    {y_x}  }
\sqrt{  1 +     y_x  ^2    }      }
\\
\\
=
{ \displaystyle   g
\left(
 { \dot{y}    -  {y_x}
\over
1 +    \dot{y}   {y_x}     } ,
 { \dot{y}   -  \dot{y}  _-
\over
1 +    \dot{y}   \dot{y}  _-        }
\right)  }
\\
\end{array}
 \\
  &  &  &  &  \\
\cline{2-5}
  &  &  &  &  \\
&
\mbox{\bf A}_{3,7}  ^b   (22)
&
\begin{array}{l}
{ \displaystyle
X_1 = { \dy} ,
\
X_2 = x { \dy} ;  }
\\
{ \displaystyle
X_3 =   ( 1   +  x ^2 )  {\dx}   +  ( x + b)   y {\dy}  ,  }
\\
{ \displaystyle
b \geq 0  }
\\
\end{array}
&
 \ddot{y}  =    { \displaystyle    { e^{ b   \arctan  x  }   \over   ( 1 + x^2 ) ^{3/2}  }  }
{ \displaystyle   f   (   u_7  , v_7 )  }
 &
{ \displaystyle   x_- = { x -     g   (   u_7  , v_7 )    \over 1 + x  g  (   u_7  , v_7 )      }  }
 \\
  &  &  &  &  \\
\hline
  &  &  &  &  \\
{\sso {sl}}  (2, \mathbb{R})
&
\mbox{\bf A}_{3,8}  (18)
&
\begin{array}{l}
  {\displaystyle
X_1 = {\dy} ,
\
X_2 =  x {\dx} +  y   {\dy} ,
 }
\\
 {\displaystyle
X_3 =   2 x y  {\dx}  +   y ^2  {\dy}     }
\\
\end{array}
&
{ \displaystyle
  \ddot{y}
=   -  {  \dot{y}  \over 2x  }
+   { \dot{y} ^3  \over x  }
f \left(
{ 1 \over  \dot{y}  }   -   {  2  x  \over  \Delta y   }  ,
   { 1 \over  \dot{y} _-  }   +   {  2  x _-   \over  \Delta y   }
\right)
}
&
{ \displaystyle
    x_-     =  {   (\Delta y)^2    \over x }
g \left(
{ 1 \over  \dot{y}  }   -   {  2  x  \over  \Delta y   }  ,
   { 1 \over  \dot{y} _-  }   +   {  2  x _-   \over  \Delta y   }
\right)
}
\\
  &  &  &  &  \\
\cline{2-5}
  &  &  &  &  \\
&
\mbox{\bf A}_{3,9}  (2)
&
\begin{array}{l}
   {\displaystyle
X_1 = {\dy}  ,
\
X_2 =  x  {\dx}  +  y  {\dy}  ,
  }
 \\
  { \displaystyle
X_3 =   2 x y  {\dx}  +   ( y ^2  - x^2 )  {\dy}  }
\\
\end{array}
&
{ \displaystyle        \ddot{y}
=     { ( 1 + \dot{y} ^2)    \dot{y} \over  x  }
+ {    ( 1 + \dot{y} ^2) ^{3/2}   \over  x  }  f ( u_9 ,  v_9 )
 }
&
{   ( x  - x_- )^2   +   ( \Delta y )^2     }       =       x x_-   g  ( u_9 ,  v_9 )
\\
  &  &  &  &  \\
\cline{2-5}
  &  &   &  &  \\
&
\mbox{\bf A}_{3,10}  (17)
&
\begin{array}{l}
 {\displaystyle
X_1 = {\dy} ,
\
X_2 =  x  {\dx}  +  y  {\dy} ,
  }
\\
 {\displaystyle
X_3 =   2 x y  {\dx}  +   ( y ^2  + x^2 )   {\dy}    }
\\
\end{array}
&
 { \displaystyle     \ddot{y}
= { | 1 - \dot{y} ^2 |     \dot{y} \over  x  }
+  {    | 1 - \dot{y} ^2| ^{3/2}   \over  x  }  f  ( u_{10} ,  v_{10}  )
 }
&
{     ( x  - x_- )^2    -  ( \Delta y )^2    }       =       x x_-   g   (  u_{10} ,  v_{10} )
\\
  &  &  &  &  \\
\cline{2-5}
  &  &  &  &  \\
&
\mbox{\bf A}_{3,11}  (11)
&
 {\displaystyle
X_1 =   {\dy} ,
\
X_2 = y  {\dy} ,
\
X_3 = y^2  {\dy}   }
&
{ \displaystyle
  \ddot{y}
= 2    { \dot{y} ^2  \over  \Delta y }
+  \dot{y}  f  \left(  x ,   { ( \Delta y )  ^2 \over  \dot{y} \dot{y}_- }   \right)
}
&
{ \displaystyle
x_- =    g \left(  x ,   { ( \Delta y )  ^2 \over  \dot{y} \dot{y}_- }   \right)
}
\\
  &  &  &  &  \\
\hline
\end{array}
$
\end{sideways}

\eject

\begin{sideways}
$
\begin{array}{|c|c|l|c|c|}
\hline
  &  &  &  &  \\
\mbox{Lie algebra} & \mbox{Case} & \mbox{Operators}  & \mbox{DODE} & \mbox{Delay relation}  \\
  &  &  &  &  \\
\hline
  &  &  &  &  \\
{\sso {o}}  (3, \mathbb{R})
&
\mbox{\bf A}_{3,12}   (3)
&
\begin{array}{l}
 {\displaystyle
X_1 = ( 1 + x^2  )  {\dx} + xy  {\dy} ,
% \quad
% X_2 = xy  {\dx} + ( 1 + y^2  )  {\dy} ,
% \quad
% X_3 = y  {\dx}  -  x  {\dy}
}
\\
 {\displaystyle
X_2 = x y  {\dx} + ( 1 + y^2  )  {\dy} ,
% \quad
% X_3 = y  {\dx}  -  x  {\dy}
}
\\
 {\displaystyle
X_3 = y  {\dx}  -  x  {\dy}   }
\\
\end{array}
&
{ \displaystyle
 \ddot{y}
=  \left( {  1 +  \dot{y} ^2 +  ( y - x  \dot{y} ) ^2 } \over  1 + x  ^2 +    y  ^2  \right)  ^{3/2}
f  ( u_{12} ,  v_{12}  )
}
&
{ \displaystyle
{
 ( x - x _- )^2
\left( 1  +     y_x   ^2
+ \left( y - x {y_x}   \right)^2   \right)
\over
(  1 + x^2 +    y ^2  )
(  1 + x _- ^2 +    y _- ^2  )
}
= g   ( u_{12} ,  v_{12}  )
}
\\
  &  &  &  &  \\
\hline
  &  &  &  &  \\
{\sso {n}}_{1,1}    \oplus {\sso {s}}_{2,1}
&
\mbox{\bf A}_{3,13}  (23)
&
\begin{array}{l}
{ \displaystyle
\left\{ X_1 = {\dx} \right\} ,
% \
% \left\{
% X_2 = {\dy} ;
% \quad
% X_3 =  y {\dy}
% \right.
}
\\
{ \displaystyle
\left\{
 X_2 = {\dy} ;
\    %\quad
X_3 =  y {\dy}
\right\}  }
\\
\end{array}
&
{ \displaystyle
\ddot{y}
=  \dot{y}
  f \left(  {   \dot{y} \over  \Delta y  } ,    {   \dot{y}_-  \over  \Delta y  }  \right)
}
&
{ \displaystyle
\Delta x  =  g \left(    {   \dot{y} \over  \Delta y  } ,    {   \dot{y}_-  \over  \Delta y  }  \right)
}
\\
   &  &  &  &  \\
\cline{2-5}
  &  &  &  &  \\
&
\mbox{\bf A}_{3,14}   (22)
&
\begin{array}{l}
{ \displaystyle
\left\{  X_1 = x {\dy}  \right\}   ,
% \
% \left\{
% X_2 = {\dy} ;
% \quad
% X_3 =   x {\dx}  + y {\dy}
% \right.
   }
\\
{ \displaystyle
\left\{
X_2 = {\dy} ;
\    %\quad
X_3 =   x {\dx}  + y {\dy}
\right\}
}
\\
\end{array}
&
{ \displaystyle
 \ddot{y}
=  {1\over x }
f \left(
\dot{y} -     {y_x}  ,
 \dot{y} -  \dot{y} _-
 \right)
 }
 &
x_-
=  x
g \left(
\dot{y} -     {y_x}  ,
 \dot{y} -  \dot{y} _-
 \right)
\\
  &  &  &  &  \\
\hline
  &  &  &  &  \\
3{\sso {n}}_{1,1}
&
\mbox{\bf A}_{3,15}  (20)
&
\begin{array}{l}
{\displaystyle
\left\{  X_1 =   {\dy}   \right\}  ,
\
\left\{  X_2 = x   {\dy}   \right\} ,
 }
\\
{\displaystyle
\left\{    X_3 =  \chi (x)   {\dy}    \right\}  ,
\
\ddot{\chi}   (x) {\not\equiv}  0  }
\\
\end{array}
&
\begin{array}{c}
 { \displaystyle
\ddot{y}
=
  {   \ddot{\chi}  (x)
\over
  \dot{\chi} (x) -   {{\chi}_x}       }   (  \dot{y} - {y_x}     )
} \\
\\
 { \displaystyle
+ \ddot{\chi}  (x)
 f   \left(
 x   ,
 { { \displaystyle    \dot{y} - {y_x}  }
\over
   \dot{\chi} (x) -   {{\chi}_x}      }
-
{ { \displaystyle          \dot{y}   -    \dot{y}_-    }
\over
   \dot{\chi} (x)     -  \dot{\chi} (x_-)     }
 \right)
 }
\end{array}
&
{ \displaystyle
 x_-
=
g
\left(
x   ,
 { { \displaystyle    \dot{y} - {y_x}  }
\over
   \dot{\chi} (x) -   {{\chi}_x}      }
- { { \displaystyle          \dot{y}   -    \dot{y}_-    }
\over
   \dot{\chi} (x)     -  \dot{\chi} (x_-)     }
\right)
}
\\
  &  &  &  &  \\
\hline
\end{array}
$
\end{sideways}

\eject

In Table~2 there were used the following notations
\begin{equation*}
u_7 =
 { \sqrt{ 1 + x^2 }     e^{ - b   \arctan  x  }    }    \left(    \dot{y}   -   {y_x}     \right) ,
\qquad
v_7 =
{ \sqrt{ 1 + x _- ^2 }     e^{ - b   \arctan (x _- ) } }   \left(    \dot{y} _-   -  {y_x}    \right)   ;
\end{equation*}
\begin{equation*}
u_9 =
 {    2 \Delta y  x   \dot{y}    -       ( \Delta y )^2  -  x_- ^2  +  x ^2
\over
 2 \Delta y  x   +      (    ( \Delta y )^2  + x_- ^2  - x ^2    )      \dot{y}  }  ,
\qquad
v_9 =
{     2 \Delta y  x_-   \dot{y}_-      +       ( \Delta y )^2  +   x ^2  -  x _- ^2
\over
 - 2 \Delta y  x_-   +      (    ( \Delta y )^2  + x ^2  - x_- ^2    )      \dot{y}_-  }    ;
\end{equation*}
\begin{equation*}
u_{10}  =
 {   2 \Delta y x  \dot{y}    -      ( \Delta y  )^2  - x^2  +  x_- ^2
\over
 2 \Delta y x    -        (      ( \Delta y  )^2  + x^2 - x_- ^2    )     \dot{y}    } ,
\qquad
v_{10}  =
 {  2 \Delta y x _- \dot{y} _-    +     ( \Delta y  )^2  + x _- ^2  -   x ^2
\over
 2 \Delta y x _-    +        (      ( \Delta y  )^2  + x _- ^2 - x ^2    )     \dot{y}  _-    }       ;
\end{equation*}
\begin{equation*}
u_{12}
 = {
 ( x - x _- )
\left( \dot{y} -  {y_x} \right)
\over
\sqrt{  1 +  \dot{y} ^2 +  ( y - x  \dot{y} ) ^2  }
\sqrt{  1 + x _- ^2 +    y _- ^2  }
}  ,
\qquad
v_{12}
= {
 ( x - x _- )
\left( \dot{y} _-  -  {y_x} \right)
\over
\sqrt{  1 +  \dot{y} _-  ^2 +  ( y _-  - x  _-  \dot{y} _-  ) ^2  }
\sqrt{  1 + x  ^2 +    y  ^2  }
}    .
\end{equation*}

\eject

$ $

\vspace{-35mm}

\subsection*{Table~3. Classification of nonlinear invariant DODSs. Symmetry algebra dimension
4.}

\begin{sideways}
$
\begin{array}{|c|c|l|c|c|}
\hline
    &   &   &   &  \\
\mbox{Lie algebra} & \mbox{Case} & \mbox{Operators} & \mbox{DODE} & \mbox{Delay relation}  \\
    &   &   &   &  \\
\hline
    &   &   &   &  \\
{\sso {n}}_{4,1}
&
\mbox{\bf A}_{4,1} (22)
&
\begin{array}{l}
  {\displaystyle
X_1 =   {\dy}  ,
\
X_2 =  x  {\dy}  ;  }
\\
 {\displaystyle
X_3 =  x ^2    {\dy}  ,
\
X_4 = {\dx}
 }
\\
\end{array}
&
{ \displaystyle
 \ddot{y}
=  { \dot{y} - \dot{y}_- \over \Delta x }
+ f \left(
 \dot{y}  +    \dot{y} _-   -   2  {y_x}
\right)
}
&
{ \displaystyle
  { \Delta x  }
= g \left(
 \dot{y}  +    \dot{y} _-   -   2  {y_x}
\right)
}
\\
    &   &   &   &  \\
\hline
    &   &   &   &  \\
{\sso {s}}_{4,1}
&
\mbox{\bf A}_{4,2}  (22)
&
\begin{array}{l}
  {\displaystyle
X_1 =    {\dy}  ,
\
X_2 =  x   {\dy}  ,
\
X_3 =  e ^{x}   {\dy}  ; }
\\
 {\displaystyle
X_4 =  {\dx}
%  ,  \quad   a \neq 0
}
\\
\end{array}
&
{ \displaystyle
\ddot{y}
=  {  e ^{ \Delta x }  \over  e ^{ \Delta x } - 1  } ( \dot{y} - \dot{y}_- )
+  f (z_{2})
}
&
{ \displaystyle
 { \Delta x  }
= g (z_{2})
}
\\
    &   &   &   &  \\
\hline
    &   &   &   &  \\
{\sso {s}}_{4,2}
&
\mbox{\bf A}_{4,3}  (22)
&
\begin{array}{l}
  {\displaystyle
X_1 =   {\dy}  ,
\
X_2 =  x {\dy}  ,
\
X_3 =  x ^2    {\dy}  ; }
\\
 {\displaystyle
X_4 = {\dx}  + y {\dy}
 }
\\
\end{array}
&
{ \displaystyle
 \ddot{y}
=  { \dot{y} - \dot{y}_- \over \Delta x }
+   e  ^{x}  f \left(
e  ^{ -x} \left(
\dot{y}  +    \dot{y} _-  -   2   {y_x}
\right)  \right)
}
&
{ \displaystyle
 { \Delta x  }
= g \left(
e  ^{ -x} \left(
\dot{y}  +    \dot{y} _-  -   2  {y_x}
\right)  \right)
}
\\
    &   &   &   &  \\
\hline
    &   &   &   &  \\
{\sso {s}}_{4,3}
&
\mbox{\bf A}_{4,4}  ^{a, \alpha}   (22)
&
\begin{array}{l}
  {\displaystyle
X_1 =    {\dy}   ,
\
X_2 =  x {\dy}  ,
\
X_3 =   |x|  ^{ \alpha }   {\dy}  ;
}
\\
  {\displaystyle
X_4 = ( 1 - a) x   {\dx}   + y   {\dy}  ,
}
\\
 {\displaystyle
a \in [-1 , 0) \cup  (0,1)  ,
\
 \alpha  \neq  \{  0, {1 \over 1-a} , 1 \}   }
\\
{\displaystyle
\mbox{(see~\cite{SWbook} for additional restriction on  $a$ and $\alpha$)} }
\\
\end{array}
&
\begin{array}{c}
{ \displaystyle
 \ddot{y}
=
 { ( \alpha -1 )  |x| ^{ \alpha -2}   ( \dot{y} - \dot{y}_- )
\over
  |x |^{ \alpha-1} \mbox{sgn} ( x )  -  |x_- | ^{ \alpha-1}  \mbox{sgn} ( x_-)  }
}
\\
\\
{ \displaystyle
+   x ^{  {2a - 1 \over 1 - a }  }       f (z_{4})
}
\\
\end{array}
&
{ \displaystyle
  x_-
= x g (z_{4})
}  \\
%   &  &   \\
% &  &
% \mbox{\bf OLD version:} \quad  {\displaystyle
% a \neq \{0,1 \},  \     \alpha  \neq  \{  0, {1 \over 1-a} \}  }
% \\
    &   &   &   &  \\
% \cline{2-5}
%     &   &   &   &  \\
% &
% \mbox{\bf A}_{4,5} (21)
% &
% \begin{array}{l}
%   {\displaystyle
% X_1 = {\dy}    ,
% \
% X_2 =  x  {\dy}  ,
% \
% X_3 =  \chi (x)   {\dy}  ; }
% \\
%   {\displaystyle
% X_4 = y   {\dy}   ,
% \
% \ddot{\chi}(x) {\not\equiv} 0   }
% \\
% \end{array}
% &
% \mbox{\bf Linear (to remove this case  ???) }
% &
% x_- = g(x)
% \\
%     &   &   &   &  \\
\hline
    &   &   &   &  \\
{\sso {s}}_{4,4}
&
\mbox{\bf A}_{4,6} ^a (22)
&
\begin{array}{l}
  {\displaystyle
X_1 =    {\dy}  ,
\
X_2 =  x    {\dy}   ,
\
X_3 =  e ^{ a x}   {\dy}   ; }
\\
{\displaystyle
X_4 =  {\dx}   +  y {\dy}  ,
\
a \neq 0 , 1 }
\\
\end{array}
&
{ \displaystyle
\ddot{y}
=  {  a   e ^{ a \Delta x }  \over  e ^{ a \Delta x } - 1  } ( \dot{y} - \dot{y}_- )
+  e^{x}   f (z_{6} )
}
&
{ \displaystyle
{ \Delta x  }
= g(z_{6} )
}
\\
    &   &   &   &  \\
\hline
    &   &   &   &  \\
{\sso {s}}_{4,5}
&
\mbox{\bf A}_{4,7}  ^ {\alpha, \beta}  (22)
&
\begin{array}{l}
  {\displaystyle
X_1 =  {\dy}   ,   }
\\
 {\displaystyle
X_2 =  e ^{ \alpha  x} \cos ( \beta x )  {\dy}  ,
\
X_3 =  e ^{ \alpha x} \sin ( \beta x )   {\dy}  ;  }
\\
{\displaystyle
X_4 =  {\dx}    + y  {\dy}  ,   }
\\
  {\displaystyle
  %\quad
 \beta \neq 0  }
\\
\end{array}
&
\begin{array}{c}
{ \displaystyle
 \ddot{y}
=  \left( { \beta \cos ( \beta \Delta x ) \over \sin ( \beta \Delta x ) }   +   \alpha  \right)  \dot{y}
-   { \beta  e ^{ \alpha  \Delta x }  \over \sin ( \beta \Delta x ) } \dot{y}_-
}
\\
\\
+  e ^{x}   f(z_{7} )
\\
\end{array}
&
{ \displaystyle
  { \Delta x  }
= g(z_{7} )
}
\\
    &   &   &   &  \\
\hline
\end{array}
$
\end{sideways}

\eject

$  $

\vspace{-15mm}

\begin{sideways}
$
\begin{array}{|c|c|l|c|c|}
\hline
    &   &   &   &  \\
\mbox{Lie algebra} & \mbox{Case} & \mbox{Operators} & \mbox{DODE} & \mbox{Delay relation}  \\
    &   &   &   &  \\
\hline
    &   &   &   &  \\
{\sso {s}}_{4,6}
&
\mbox{\bf A}_{4,8}  (24)
&
\begin{array}{l}
  {\displaystyle
X_1 =  {\dy}   ,
\   %\quad
X_2 = {\dx}    ,
\   %\quad
X_3 = x  {\dy}  ;   }
\\
 {\displaystyle
  %\quad
X_4 =  x   {\dx}
  }
\\
\end{array}
  &
{ \displaystyle
   \ddot{y}
= { 1 \over (  \Delta x ) ^2  } f \left(
 \Delta x
\left(  \dot{y}_-  -  {y_x}  \right)
\right)
}
 &
{ \displaystyle
\Delta x
= {
g \left(
 \Delta x
\left(  \dot{y}_-  -  {y_x}  \right)
\right)
\over
 \dot{y} -  {y_x}
}
}
\\
    &   &   &   &  \\
\hline
    &   &   &   &  \\
{\sso {s}}_{4,8}
&
\mbox{\bf A}_{4,9} ^a  (24)
&
\begin{array}{l}
  {\displaystyle
X_1 =  {\dy}   ,
\   %\quad
X_2 =  {\dx}   ,
\   %\quad
X_3 = x  {\dy}  ;  }
\\
 {\displaystyle
   %\quad
X_4 =  x  {\dx}   +  a y {\dy}  ,
\   %\quad
a \neq  0,1
  }
\\
\end{array}
&
 { \displaystyle
  \ddot{y}
=  |  \Delta x | ^ { \alpha -2  }   f \left(
 |  \Delta x | ^ { 1- \alpha   }
\left(  \dot{y}_-  -  {y_x}  \right)
\right)
}
&
{ \displaystyle
|  \Delta x | ^ { 1 -  \alpha  }
=
{
g \left(
 |  \Delta x | ^ { 1- \alpha   }
\left(  \dot{y}_-  -  {y_x}  \right)
\right)
\over
  \dot{y} -  {y_x}
}
}
\\
    &   &   &   &  \\
\hline
    &   &   &   &  \\
{\sso {s}}_{4,10}
&
\mbox{\bf A}_{4,10}  (25)
&
\begin{array}{l}
  {\displaystyle
X_1 =  {\dy}   ,
\   %\quad
X_2 =  {\dx}    ,
\   %\quad
X_3 = x  {\dy}  ; }
 \\
  {\displaystyle
   %\quad
X_4 =    x   {\dx}  +    ( 2 y + x^2 )  {\dy}
}
\\
\end{array}
&
  { \displaystyle
 \ddot{y}
=  {  \dot{y}  -  \dot{y}_-  \over     \Delta x   }
+  f \left(
( \Delta x )    \mbox{exp}   \left[
{ \displaystyle     \dot{y}_-  -    y_x   \over  \Delta x }
\right]
\right)
}
&
{ \displaystyle
{  \Delta x   \over \mbox{exp}   \left[
{ \displaystyle  {   \dot{y} -  {y_x}   \over  \Delta x }  }
\right]  }
=  g \left(
( \Delta x )     \mbox{exp}   \left[
{ \displaystyle  \dot{y}_-  -  {y_x}    \over   \Delta x }
\right]
\right)
}
\\
    &   &   &   &   \\
\hline
    &   &   &   &  \\
{\sso {s}}_{4,11}
&
\mbox{\bf A}_{4,11} (24)
&
\begin{array}{l}
  {\displaystyle
X_1 =  {\dy}   ,
\   %\quad
X_2 = {\dx}   ,
\   %\quad
X_3 = x  {\dy}  ;  }
\\
 {\displaystyle
  %\quad
X_4 =  x   {\dx}     +  y  {\dy}
  }
\\
\end{array}
&
{ \displaystyle
     \ddot{y}
=  { 1\over  \Delta x  } f \left(
\dot{y}_-  -  {y_x}
\right)
}
&
{ \displaystyle
    \dot{y} -  {y_x}
=  g \left(
\dot{y}_-  -  {y_x}
\right)
}
\\
    &   &   &   &  \\
\cline{2-5}
    &   &   &   &  \\
&
\mbox{\bf A}_{4,12} (23)
&
\begin{array}{l}
  {\displaystyle
X_1 =  {\dy}   ,
\   %\quad
X_2 =  x  {\dy}  ,
\   %\quad
X_3 =   {\dx}    ;  }
\\
 {\displaystyle
  %\quad
X_4 =  y {\dy}
 }
\\
\end{array}
&
{ \displaystyle
\ddot{y}
= \left( \dot{y} - {y_x} \right)
f \left(
{ \displaystyle \dot{y}_-   - {y_x} \over \displaystyle   \dot{y} - {y_x}  }
\right)
}
&
{ \displaystyle
 \Delta  x
= g \left(
{ \displaystyle \dot{y}_- - {y_x} \over \displaystyle   \dot{y} - {y_x}  }
\right)
}
\\
    &   &   &   &  \\
\hline
    &   &   &   &  \\
{\sso {s}}_{4,12}
&
\mbox{\bf A}_{4,13}  (4)
&
\begin{array}{l}
  {\displaystyle
X_1 =   {\dx}  ,
\   %\quad
X_2 =   {\dy}  ;
 }
\\
 {\displaystyle
X_3 = x  {\dx}  + y {\dy} ,
\
X_4 = y  {\dx}  -  x  {\dy}
}
\\
\end{array}
&
{ \displaystyle
     \ddot{y}
=   { ( 1 +  \dot{y} ^2 ) ^2  \over   \Delta x \left(   \dot{y}   -  {  \displaystyle {y_x}   }  \right) }
f \left(
 {  \dot{y}  -   \dot{y}_-     \over 1 +  \dot{y}   \dot{y}_-   }
\right)
}
&
{ \displaystyle
{   \dot{y}  -  {  \displaystyle  {y_x} }    \over
   1    + \dot{y}  {  \displaystyle  {y_x}  }   }
= g \left(
 {  \dot{y}  -   \dot{y}_-     \over 1 +  \dot{y}   \dot{y}_-   }
\right)
}
\\
    &   &   &   &  \\
\cline{2-5}
    &   &   &   &  \\
&
\mbox{\bf A}_{4,14}   (23)
&
\begin{array}{l}
  {\displaystyle
X_1 =  {\dy}  ,
\   %\quad
X_2 =  x  {\dy} ;
 }
\\
 {\displaystyle
X_3 =  y {\dy} ,
\
X_4 = ( 1 + x^2 )   {\dx}  +   x y {\dy}
  }
\\
\end{array}
&
{ \displaystyle
 \ddot{y}
=   { \dot{y} -  {  \displaystyle {y_x}    }  \over  ( x ^2 + 1 ) }
f  (  z )  ,
%  \sqrt{ { x _- ^2 + 1 \over x  ^2 + 1 } }
% { \displaystyle \dot{y}_- - {y_x} \over  \displaystyle  \dot{y} - {y_x}  }
% \right)
\
z
=
 \sqrt{ { x _- ^2 + 1 \over x  ^2 + 1 } }
{ \displaystyle \dot{y}_- - {y_x} \over  \displaystyle  \dot{y} - {y_x}  }
}
&
{ \displaystyle
x_- =  {   x -   g (z)   \over 1 + x g (z)  }
}
\\
   &  &  &  &  \\
\hline
\end{array}
$
\end{sideways}

\eject

\begin{sideways}
$
\begin{array}{|c|c|l|c|c|}
\hline
    &   &   &   &  \\
\mbox{Lie algebra} & \mbox{Case} & \mbox{Operators} & \mbox{DODE} & \mbox{Delay relation}  \\
    &   &   &   &  \\
\hline
    &   &   &   &  \\
{\sso {n}}_{1,1} \oplus {\sso {s}}_{3,1}
&
\mbox{\bf A}_{4,15} ^a  (22)
&
\begin{array}{l}
  {\displaystyle
\left\{    X_1 =   |x|  ^{1 \over 1-a}   {\dy}   \right\}  ,  }
\\
 {\displaystyle
   %\quad
\left\{
X_2 =   {\dy}  ,
\   %\quad
X_3 =  x {\dy}  ;
\
X_4 = ( 1 - a ) x   {\dx}  + y   {\dy}
\right\}
}  ,
\\
{\displaystyle
a \in [-1 , 0) \cup  (0,1)   }
% OLD version:     a \neq \{0,1 \}
%   \alpha  = {1 \over 1-a}  }
\\
\end{array}
&
\begin{array}{c}
{ \displaystyle
 \ddot{y}
=  {  {  a  \over 1 - a }  |x| ^{  2a -1  \over 1 - a }      ( \dot{y} - \dot{y}_- )
\over
    |x |^{  a  \over 1 - a }   \mbox{sgn} ( x )
-  |x_- | ^{  a  \over 1 - a }   \mbox{sgn} ( x_-)  }
}
 \\
\\
{ \displaystyle
+  x ^{ 2a - 1 \over 1 - a }      f (z_{15})   }
\\
\end{array}
&
{ \displaystyle
 x_-  = x g (z_{15})
}
\\
   &   &   &   &  \\
\hline
    &   &   &   &  \\
{\sso {n}}_{1,1} \oplus {\sso {s}}_{3,2}
&
\mbox{\bf A}_{4,16} (22)
&
\begin{array}{l}
  {\displaystyle
\left\{   X_1 =  e ^{x}   {\dy}   \right\} , }
\\
 {\displaystyle
   %\quad
\left\{
X_2 =    {\dy}   ,
\   %\quad
X_3 =  x    {\dy}   ;
\   %\quad
X_4 =   {\dx}  + y {\dy}
\right\}
  }
\\
\end{array}
&
{ \displaystyle
 \ddot{y}
=  {     e ^{  \Delta x }  \over  e ^{  \Delta x } - 1  } ( \dot{y} - \dot{y}_- )
+  e^{x} f (z_{16} )
}
&
{ \displaystyle
 { \Delta x  } = g (z_{16}  )
}
\\
    &   &   &   &  \\
\hline
    &   &   &   &  \\
{\sso {n}}_{1,1} \oplus {\sso {s}}_{3,3}
&
\mbox{\bf A}_{4,17} ^{ \alpha } (22)
&
\begin{array}{l}
  {\displaystyle
\left\{  X_1 = {\dy}  \right\}   ,  }
\\
 {\displaystyle
   %\quad
\left\{
X_2 =    {\dx}   ,
\   %\quad
X_3 =  e ^{ \alpha  x} \cos  x   {\dy}  ,
\
X_4 =  e ^{ \alpha x}  \sin   x  {\dy}
\right\} ,
}
\\
 {\displaystyle
   %\quad
{ \alpha  \geq 0 }
}
\\
\end{array}
&
{ \displaystyle
 \ddot{y}
=  \left( {  \cos ( \Delta x ) \over \sin ( \Delta x ) }   +   \alpha  \right)  \dot{y}
- { e ^{ \alpha  \Delta x }  \over \sin (  \Delta x ) } \dot{y}_-
+  f(z_{17} )
}
&
{ \displaystyle
 { \Delta x  } = g (z_{17} )
}
\\
    &   &   &   &  \\
\hline
\end{array}
$
\end{sideways}

\eject

\begin{sideways}
$
\begin{array}{|c|c|l|c|c|}
\hline
    &   &   &   &  \\
\mbox{Lie algebra} & \mbox{Case} & \mbox{Operators} & \mbox{DODE} & \mbox{Delay relation}  \\
    &   &   &   &  \\
\hline
    &   &   &   &  \\
{\sso {n}}_{1,1} \oplus {\sso {sl}}  (2, \mathbb{R})
&
\mbox{\bf A}_{4,18} (14)
&
\begin{array}{l}
  {\displaystyle
\left\{     X_1 = {\dx}      \right\}   , }
\\
  {\displaystyle
   %\quad
\left\{
X_2 =  {\dy}  ,
\   %\quad
X_3 = y {\dy}  ,
\   %\quad
X_4  = y^2 {\dy}
\right\}  }
\\
\end{array}
&
{ \displaystyle
 \ddot{y}
=   2 {  \dot{y} ^2   \over  \Delta  y  }
+   \dot{y} f \left(
{  (  \Delta  y ) ^2 \over \dot{y}  \dot{y} _- }
\right)
}
&
{ \displaystyle
\Delta  x
= g \left(
{  (  \Delta  y ) ^2 \over \dot{y}  \dot{y} _- }
\right)
}
\\
    &   &   &   &  \\
\cline{2-5}
    &   &   &   &  \\
&
\mbox{\bf A}_{4,19}  (19)
&
\begin{array}{l}
  {\displaystyle
\left\{ X_1 =   x   {\dx}      \right\} , }
\\
 {\displaystyle
   %\quad
\left\{
X_2 = {\dy}  ,
\   %\quad
X_3 = x   {\dx}   +   y {\dy} , \right. }
\\
{\displaystyle
  %\quad
\left.
X_4 = 2 x y   {\dx}   +  y^2   {\dy}
\right\}  }
\\
\end{array}
&
\begin{array}{c}
{ \displaystyle
     \ddot{y}
=
 -   {    \dot{y}  \over   2 x   }
+ { x  _-   \dot{y} ^3   \over  ( \Delta y )^2 }
}
\\
\\
{ \displaystyle
\times  f \left(
{ ( \Delta y )^2  \over x   x  _- }
\left(   { 1 \over  \dot{y}_-  }  +  { 2 x _- \over  \Delta y } \right) ^2
\right)
}
\\
\end{array}
&
\begin{array}{c}
{ \displaystyle
 { ( \Delta y )^2  \over x   x  _- }
\left(   { 1 \over  \dot{y} }  -  { 2 x  \over  \Delta y } \right) ^2
}
\\
\\
 { \displaystyle
 = g  \left(
{ ( \Delta y )^2  \over x   x  _- }
\left(   { 1 \over  \dot{y}_-  }  +  { 2 x _- \over  \Delta y } \right) ^2
\right)
}
\\
\end{array}
\\
    &   &   &   &  \\
\hline
    &   &   &   &  \\
2 {\sso {s}}_{2,1}            % \oplus {\sso {s}}_{2,1}
&
\mbox{\bf A}_{4,20} (13)
&
\begin{array}{l}
  {\displaystyle
\left\{
X_1 =   {\dx}   ;
\   %\quad
X_2 = x   {\dx}
\right\} ,  }
\\
{\displaystyle
   %\quad
\left\{
X_3 =  {\dy}   ;
\   %\quad
X_4  = y  {\dy}
\right\} }
\\
\end{array}
&
{ \displaystyle
 \ddot{y}
=  { \dot{y} \over   \Delta  x  }
 f  \left(
 {   \dot{y} _- \over   \dot{y}  }
\right)
}
 &
{ \displaystyle
   \Delta  x
=  {  \Delta  y \over  \dot{y}  }
g  \left(
 {   \dot{y} _- \over   \dot{y}  }
\right)
}
\\
    &   &   &   &  \\
\cline{2-5}
    &   &   &   &  \\
&
\mbox{\bf A}_{4,21}  (23)
&
\begin{array}{l}
  {\displaystyle
\left\{
X_1 =  {\dy}  ;
\   %\quad
X_2 = x   {\dx}  +   y  {\dy}
\right\} ,   }
\\
 {\displaystyle
   %\quad
\left\{
X_3 =  x   {\dy}  ;
\   %\quad
X_4 =  x   {\dx}
\right\} }
\\
\end{array}
&
{ \displaystyle
 \ddot{y}
= { \displaystyle  \dot{y} - {y_x}   \over x }
 f  \left(
 { \displaystyle \dot{y}_-   - {y_x}
\over
\displaystyle \dot{y} - {y_x}  }
\right)
}
&
{ \displaystyle
 x_-
= x g  \left(
 { \displaystyle \dot{y}_-   - {y_x}
\over
\displaystyle \dot{y} - {y_x}  }
\right)
}
\\
    &   &   &   &  \\
\hline
%     &   &   &   &  \\
% 4 {\sso {n}}_{1,1}
% &
% \mbox{\bf A}_{4,22}  (20)
% &
%   {\displaystyle
% \left\{  X_1 = {\dy}  \right\}  ,
% \   %\quad
% \left\{  X_2 = x  {\dy} \right\} ,
% \   %\quad
% \left\{  X_3 =  \chi_1 (x)  {\dy}  \right\} ,
% \   %\quad
% \left\{  X_4 =  \chi_2 (x)  {\dy}  \right\}
% }
% &
% \mbox{Linear (see T� ???)}
% &
% x_- = g(x)
% \\
%     &   &   &   &  \\
% &
% &
% \mbox{$1$, $x$, $\chi_1 (x)$ and $\chi_2 (x) $ are linearly independent}
% &
% &
% \\
%  &  &  &  &  \\
% \hline
\end{array}
$
\end{sideways}

\eject

In Table~3 there were used the notations
\begin{equation*}
z _{2}
=  \left(  { e ^{ \Delta x } - 1 \over  \Delta x }  -   1  \right)    \dot{y}
+  \left(  e ^{ \Delta x }  - { e ^{ \Delta x } - 1 \over  \Delta x }    \right)    \dot{y} _-
-      (  e ^{ \Delta x } - 1  )   {y_x}     ,
\end{equation*}
\begin{multline*}
z _{4} = x ^{ { 2a - 1 \over a-1} - \alpha }
\left[
{ \alpha}  |x_-| ^{ \alpha -1 }   \mbox{sgn} ( x_-)
         (   \dot{y} -   {y_x}   )
\right.
\\
\left.
-
{ \alpha} | x |^{ \alpha -1 }  \mbox{sgn} ( x )
           (  \dot{y}_-   -  {y_x}  )
-
{ |x| ^{ \alpha}  - |x _-| ^{ \alpha} \over   x - x_- }    ( \dot{y}   -  \dot{y}_-  )
 \right]    ,
\end{multline*}
\begin{equation*}
z _{6}
=   e^{ -x} \left[
\left(  { e ^{ a \Delta x } - 1 \over  \Delta x } -   a  \right)    \dot{y}
+  \left(   a e ^{ a \Delta x }  - { e ^{ a \Delta x } - 1 \over  \Delta x }    \right)    \dot{y} _-
-    a (  e ^{ a \Delta x } - 1  )   {y_x}
  \right]   ,
\end{equation*}
\begin{multline*}
z_{7}
= e ^{ -x} \left[
\left(  { \beta  \cos ( \beta \Delta x )  -  \beta   e ^{ - \alpha  \Delta x }   \over  \sin ( \beta \Delta x )  }
-  \alpha   \right)  \dot{y}
\right.
\\
\left.
+
\left(  { \beta  \cos ( \beta \Delta x )  -  \beta   e ^{ \alpha  \Delta x } \over   \sin ( \beta \Delta x ) }
+  \alpha      \right)    \dot{y} _-
+    ( \alpha ^2 + \beta ^2 )    ( y - y _- )
\right]   ,
\end{multline*}
\begin{multline*}
z_{15}   = x ^{ { 2a  \over a-1} }
\left[
     |x_-| ^ {  a  \over 1 - a }      \mbox{sgn} ( x_-)   \left(   \dot{y}   -   {y_x} \right)
\right.
\\
\left.
-   | x |^ {  a  \over 1 - a }    \mbox{sgn} ( x )  \left(   \dot{y}_-  -  {y_x}  \right)
   - ( 1 - a ) { |x| ^ {  1  \over 1 - a }  - |x _-| ^{  1  \over 1 - a }  \over   x - x_- }    \left(  \dot{y}  - \dot{y}_- \right)
\right]   ,
\end{multline*}
\begin{equation*}
z _{16} =   e^{ -x} \left[
\left(   { e ^{ \Delta x } - 1 \over  \Delta x }    -   1    \right)    \dot{y}
+  \left(   e ^{  \Delta x }  - { e ^{  \Delta x } - 1 \over  \Delta x }    \right)    \dot{y} _-
-     (  e ^{ \Delta x } - 1  )   {y_x}
  \right]    ,
\end{equation*}
\begin{equation*}
z _{17}
=  \left(  { \cos ( \Delta x )  -   e ^{ - \alpha  \Delta x }   \over  \sin ( \Delta x )  }
-  \alpha   \right)  \dot{y}
+
 \left(  { \cos (  \Delta x ) -   e ^{ \alpha  \Delta x } \over   \sin ( \Delta x ) }
+  \alpha      \right)    \dot{y} _-
+    ( \alpha ^2 + 1 )    ( y - y _- )    .
\end{equation*}

 \eject

\subsection*{Table~4. Classification of nonlinear invariant DODSs.
Symmetry algebra dimensions 5 and 6. }

\begin{sideways}
$
\begin{array}{|c|c|l|c|c|}
\hline
  &  &  &  &  \\
\mbox{Lie algebra} & \mbox{Case} & \mbox{Operators}  & \mbox{DODE} & \mbox{Delay relation}  \\
  &  &  &  &  \\
\hline
  &  &  &  &  \\
{\sso {s}}_{5,33}
&
\mbox{\bf A}_{5,1} (24)
&
\begin{array}{l}
 {\displaystyle
X_1 =  {\dy}   ,
\    %quad
X_2 = x  {\dy}  ,
\    %quad
X_3 = x ^2  {\dy}   , }
\\
{\displaystyle
   %quad
X_4 =   {\dx}   ;
\    %quad
X_5 = x   {\dx}
 }
\\
\end{array}
  &
{ \displaystyle
 \ddot{y}
=   2  {      \dot{y} -  {y_x}      \over \Delta x}
+  { C_1 \over  (  \Delta x ) ^  2  }
}
 &
{ \displaystyle
 \Delta x
= {  C_2  \over  \dot{y} +  \dot{y}_-   - 2 {y_x}  }
}
\\
  &  &  &  &  \\
\hline
  &  &  &  &  \\
 {\sso {s}}_{5,34}
 &
\mbox{\bf A}_{5,2} (25)
 &
\begin{array}{l}
 {\displaystyle
X_1 =  {\dy}   ,
\    %quad
X_2 = x  {\dy}  ,
\    %quad
X_3 = x ^2 {\dy}  ,  }
\\
{\displaystyle
   %quad
X_4 =  {\dx}    ; }
\
 {\displaystyle
X_5 =    x   {\dx}  +    ( 3 y + x^3 )  {\dy}
  }
\\
\end{array}
&
{ \displaystyle
   \ddot{y}
=
 {    4  \dot{y} + 2   \dot{y}_-  -   6  {y_x}    \over  \Delta x  }
+  C_1  \Delta x
}
&
{ \displaystyle
 { 1 \over (  \Delta x ) ^ 2  }
\left(  \dot{y} +  \dot{y}_-   - 2 {y_x}  \right)
- \ln  | \Delta x |
= C_2
}
\\
  &  &  &  &  \\
\hline
  &  &  &  &  \\
{\sso {s}}_{5,35}
 &
 \mbox{\bf A}_{5,3} (24)
 &
\begin{array}{l}
 {\displaystyle
X_1 =  {\dy}   ,
\    %quad
X_2 = x  {\dy} ,
\    %quad
X_3 = x ^2  {\dy}  , }
\\
{\displaystyle
    %quad
X_4 =  {\dx}     ;
  }
\
 {\displaystyle
X_5 =x   {\dx}   +  \alpha y  {\dy}  ,
\    %quad
 \alpha \neq \{ 0,2\}
  }
\\
\end{array}
&
{ \displaystyle
 \ddot{y}
=  2  {     \dot{y} -  {y_x}     \over    \Delta x   }
+  C_1  | \Delta x | ^ { \alpha - 2 }
}
&
{ \displaystyle
 | \Delta x | ^ { \alpha -1 }
=   C_2
\left(  \dot{y} +  \dot{y}_-   - 2 {y_x}  \right)
}
\\
  &  &  &  &  \\
\hline
  &  &  &  &  \\
 {\sso {s}}_{5,36}
 &
\mbox{\bf A}_{5,4} (24)
 &
\begin{array}{l}
 {\displaystyle
X_1 =  {\dy} ,
\    %quad
X_2 = x  {\dy}  ,
\    %quad
X_3 = x ^2  {\dy}   , }
\\
 {\displaystyle
    %quad
X_4 = {\dx} ;
  }
\
 {\displaystyle
X_5 =x  {\dx}  + 2 y  {\dy}
  }
\\
\end{array}
&
{ \displaystyle
 \ddot{y}
=
2  {  \dot{y} -  {y_x}   \over   \Delta x  }
+  C_1
}
&
{ \displaystyle
  \Delta x  = C_2 \left(  \dot{y} +  \dot{y}_-   - 2 {y_x}  \right)
}
\\
  &  &  &  &  \\
\hline
%   &  &  &  &  \\
%  {\sso {s}}_{5,37}
%   &
%  \mbox{\bf A}_{5,5} (24)
%  &
%  {\displaystyle
% X_1 =  {\dy}   ,
% \    %quad
% X_2 = x  {\dy}  ,
% \    %quad
% X_3 = x ^2 {\dy}   ,
% \    %quad
% X_4 = {\dx}   ;
% \    %quad
% X_5 = y  {\dy}
% }
% &
% { \displaystyle
%   \ddot{y}
% = 2  {     \dot{y} -  {y_x}   \over     \Delta x   }
% +  C_1 \left(  \dot{y} +  \dot{y}_-   - 2 {y_x} \right)
% }
% &
% { \displaystyle
% \Delta x
% = C_2
% }
% \\
%   &  &  &  &  \\
% \hline
  &  &  &  &  \\
  {\sso {s}}_{5,44}
  &
 \mbox{\bf A}_{5,5} (26)
  &
\begin{array}{l}
 {\displaystyle
X_1 =  {\dy}   ,
\    %quad
X_2 = x {\dy} ,
\    %quad
X_3 = {\dx}   ;    }
\\
{\displaystyle
    %quad
X_4 = x  {\dx}  ,
\    %quad
X_5 = y  {\dy}
 }
\\
\end{array}
&
{ \displaystyle
  \ddot{y}
= C_1  {   \dot{y} -  {y_x}    \over    \Delta x   }
}
&
{ \displaystyle
  \Delta x ={ ( 1 - C_2)  \Delta y \over  \dot{y} - C_2 \dot{y}_- }
}
\\
  &  &  &  &  \\
\hline
\end{array}
$
\end{sideways}

\eject

% $  $

% \vspace{-20mm}

\begin{sideways}
$
\begin{array}{|c|c|l|c|c|}
\hline
  &  &  &  &  \\
\mbox{Lie algebra} & \mbox{Case} & \mbox{Operators}  & \mbox{DODE} & \mbox{Delay relation}  \\
  &  &  &  &  \\
\hline
  &  &  &  &  \\
  {\sso {sl}}  (2, \mathbb{R})  \ltimes  2  {\sso {n}}_{1,1}
  &
\mbox{\bf A}_{5,6} (27)
 &
\begin{array}{l}
 {\displaystyle
\left\{
X_1 = {\dx}   ,
\    %quad
X_2 = 2 x  {\dx}+  y   {\dy}  ,
\    %quad
X_3 = x^2  {\dx} + x y   {\dy}
\right\}
  }
\\
 {\displaystyle
\ltimes
\left(
\left\{
X_4 =  {\dy}
\right\} ,
\    %quad
\left\{
X_5 = x   {\dy}
\right\}
\right)
  }
\\
\end{array}
&
{ \displaystyle
     \ddot{y}  = C_1   \displaystyle   \left( \dot{y} -  {y_x}  \right) ^3
}
&
{ \displaystyle
 \Delta x
= {  C_2   \over     \left(   \dot{y} -  {y_x}  \right)    \left(  \dot{y}_-  -  {y_x} \right)    }
}
\\
  &  &  &  &  \\
\cline{2-5}
  &  &  &  &  \\
%  {\sso {sl}}  (2, \mathbb{R})  \ltimes  2  {\sso {n}}_{1,1}
  &
\mbox{\bf A}_{5,7} (5)
  &
\begin{array}{l}
 {\displaystyle
\left\{
X_1 = x  {\dx}  - y {\dy}  ,
\    %quad
X_2 = y  {\dx}  ,
\    %quad
X_3 =   x  {\dy}
\right\}
  }
\\
 {\displaystyle
\ltimes
\left(
\left\{
X_4 = {\dx}
\right\}   ,
\    %quad
\left\{
X_5 =  {\dy}
\right\}
\right)
 }
\\
\end{array}
&
{ \displaystyle
  \ddot{y}   = C_1
 ( \Delta x ) ^3 \left( \dot{y} -  {y_x} \right) ^3
}
 &
{ \displaystyle
    ( \Delta x ) ^2
= { C_2    (    \dot{y}  -  \dot{y} _-  )  \over  \left( \dot{y} -  {y_x} \right)   \left( \dot{y}_-  -  {y_x} \right)  }
}
\\
  &  &  &  &  \\
\hline
  &  &  &  &  \\
 {\sso {s}}_{2,1} \oplus {\sso {sl}}  (2, \mathbb{R})
 &
 \mbox{\bf A}_{5,8} (15)
 &
\begin{array}{l}
 {\displaystyle
 \left\{
X_1 = {\dx}  ,
\    %quad
X_2 = x {\dx}
\right\}  ,
 }
\\
 {\displaystyle
\left\{
X_3 = {\dy} ,
\    %quad
X_4  = y {\dy} ,
\    %quad
X_5 = y^2  {\dy}
\right\}
  }
\\
\end{array}
&
{ \displaystyle
      \ddot{y}
=  2  {  \dot{y}^2   \over \Delta y  }
+   C_1 {  \dot{y}   \over  \Delta x }
}
&
{ \displaystyle
 ( \Delta  x   )  ^2
=
 C_2    {   ( \Delta  y  )  ^2    \over  \dot{y}  \dot{y}_-     }
}
 \\
  &  &  &  &  \\
\hline
\hline
  &  &  &  &  \\
 {\sso {s0}}  (3, 1 )
&
 \mbox{\bf A}_{6,1} (7)
&
\begin{array}{l}
{\displaystyle
X_1 = {\dx} ,
\   %\quad
X_2 = {\dy} ,
\   %\quad
X_3 = x {\dx}  + y  {\dy} ,
\   %\quad
X_4 =  y  {\dx} - x  {\dy} ,
}
\\
{\displaystyle
X_5 = (x^2 - y^2 )  {\dx}   +  2 x y {\dy}  ,
\   %\quad
X_6 =   2 x y {\dx}  +     (y^2 - x^2 )  {\dy }
}
\\
\end{array}
&
{\displaystyle
\ddot{y}
-
2   {    \dot{y} ^2 + 1   \over   y_x^2 + 1 }
{    \dot{y}  - y_x    \over \Delta x }
= 0
}
&
\begin{array}{c}
{\displaystyle
\arctan  ( \dot{y} )   +    \arctan  ( \dot{y} _- )
}
\\
{\displaystyle
   - 2  \arctan  (  y_x  ) = C
}
\\
\end{array}
\\
   &  &  &  &  \\
\hline
   &  &  &  &  \\
  {\sso {sl}}  (2, \mathbb{R})  \ltimes  3  {\sso {n}}_{1,1}
&
\mbox{\bf A}_{6,2} (27)
&
\begin{array}{l}
 {\displaystyle
\left\{
X_1 = {\dx}  ,
\    %quad
X_2 =  x  {\dx} +  y   {\dy}  ,
\    %quad
X_3 = x^2  {\dx} + 2 x y  {\dy}
\right\}
}
\\
 {\displaystyle
\ltimes
\left(
\left\{
X_4 = {\dy}
\right\}    ,
\    %quad
\left\{
X_5 = x   {\dy}
\right\}  ,
\    %quad
\left\{
X_6 = x ^2   {\dy}
\right\}
\right)
}
\\
\end{array}
&
{ \displaystyle
 \ddot{y}    -   2 {  \dot{y}  -  {y_x}     \over \Delta x}  = 0
}
&
{ \displaystyle
   \dot{y}  +   \dot{y}_-  -  2 {y_x}  = C
} \\
  &  &  &  &  \\
\hline
  &  &  &  &  \\
{\sso {sl}}  (2, \mathbb{R})  \oplus {\sso {sl}}  (2, \mathbb{R})
&
\mbox{\bf A}_{6,3} (16)
&
\begin{array}{l}
{\displaystyle
\left\{
X_1 = {\dx}  ,
\    %quad
X_2 = x {\dx}  ,
\    %quad
X_3 = x^2  {\dx}
\right\}    ,
}
\\
{\displaystyle
\left\{
X_4 =  {\dy}  ,
\    %quad
X_5  = y  {\dy}  ,
\    %quad
X_6 = y^2  {\dy}
\right\}
}
\\
\end{array}
&
{ \displaystyle
  \ddot{y} - 2  {  \dot{y} ^2  \over \Delta y  }  + {  2  \dot{y}   \over \Delta x }  = 0
}
&
{ \displaystyle
 (  \Delta   x )   ^2  = C {  (  \Delta   y )   ^2  \over   \dot{y}  \dot{y}_-  }
} \\
  &  &  &  &  \\
\hline
\end{array}
$
\end{sideways}

\end{document}